\documentclass{article}

\usepackage[utf8]{inputenc}
\usepackage[T1]{fontenc}

\usepackage{lmodern}
\usepackage{textcomp}

\usepackage[leqno]{amsmath}
\usepackage{amsfonts,%
   amssymb,%
   amsthm,%
   amsxtra%
 }

\usepackage{mathrsfs} 

\usepackage{mathdots}

\usepackage[all,tips,2cell]{xy}
\SelectTips{cm}{}
\UseAllTwocells

\usepackage{natbib} \setcitestyle{numbers,square}
\let\oldbibsection\bibsection
\renewcommand{\bibsection}{\oldbibsection\addcontentsline{toc}{section}{References}}
\usepackage{xspace} 
\usepackage{microtype}

\usepackage[colorlinks=true]{hyperref}

\swapnumbers
\theoremstyle{plain}
\newtheorem*{theorem*}{Theorem}
\newtheorem{theorem}[equation]{Theorem}
\newtheorem{proposition}[equation]{Proposition}
\newtheorem*{proposition*}{Proposition}
\newtheorem{lemma}[equation]{Lemma}

\newtheorem{lemma-def}[equation]{Lemma-Definition}

\theoremstyle{definition}
\newtheorem{definition}[equation]{Definition}
\newtheorem{example}[equation]{Example}
\newtheorem{assertion}[equation]{Assertion}

\newtheorem{remark}[equation]{Remark}

\numberwithin{equation}{subsection}%
\defcitealias{ButterfliesI}{Part I}

\newcommand{\eqdef}{\overset{\text{\normalfont\tiny def}}{=}}
\newcommand{\coloneq}{\mathrel{\mathop:}=}
\renewcommand{\phi}{\varphi}
\renewcommand{\epsilon}{\varepsilon}

\newcommand{\Cech}{\v{C}ech\xspace}

\newcommand{\loccit}{loc.\ cit.\xspace}

\newcommand{\ie}{i.e.\xspace}

\newcommand{\lto}{\longrightarrow} 
\newcommand{\lmto}{\longmapsto}
\newcommand{\iso}{\simeq}
\newcommand{\isoto}{\overset{\sim}{\rightarrow}} 
\newcommand{\lisoto}{\overset{\sim}{\longrightarrow}} 
\newcommand{\coin}{\equiv} 

\newcommand{\coho}[1]{\operatorname{\mathrm{#1}}}%
\renewcommand{\H}{\coho{H}} 

\DeclareMathOperator{\Aut}{Aut}
\newcommand{\shAut}{\operatorname{\underline{\mathrm{Aut}}}}

\DeclareMathOperator{\B}{B}
\DeclareMathOperator{\cech}{\check{C}}

\DeclareMathOperator{\cosk}{cosk}

\newcommand{\cprod}[1]{\operatorname{\wedge}^{#1}\displaylimits} 
\newcommand{\del}{\partial} 

\newcommand{\equ}{\operatorname{\mathscr{E}\mspace{-3mu}\mathit{q}}}
\DeclareMathOperator{\Ext}{Ext}

\DeclareMathOperator{\catExt}{\mathsf{Ext}}
\DeclareMathOperator{\shcatExt}{\mathscr{E}\mspace{-3mu}\mathit{xt}}
\DeclareMathOperator{\Hom}{Hom} 

\DeclareMathOperator{\catHom}{\mathsf{Hom}}

\DeclareMathOperator{\shHom}{\underline{\mathrm{Hom}}}
\DeclareMathOperator{\shcatHom}{\mathscr{H}\mspace{-4mu}\mathit{om}}
\DeclareMathOperator{\id}{id} 

\DeclareMathOperator{\Id}{Id}

\DeclareMathOperator{\shIsom}{\underline{\mathrm{Isom}}}

\newcommand{\colim}{\varinjlim}
\renewcommand{\lim}{\varprojlim}

\DeclareMathOperator{\W}{\overline{W}}

\newcommand{\cat}[1]{\mathsf{#1}} 
\newcommand{\bicat}[1]{\underline{\mathsf{#1}}} 
\DeclareMathOperator{\Ob}{Ob}
\DeclareMathOperator{\Mor}{Mor}

\newcommand{\Set}{\cat{Set}}

\newcommand{\site}[1]{\cat{#1}}

\newcommand{\sS}{\site{S}}
\newcommand{\s}{\sS}

\newcommand{\T}{\cat{T}}

\newcommand{\grpd}[1]{\cat{#1}}

\newcommand{\smp}[1]{\underline{#1}}

\newcommand{\stack}[1]{\mathscr{#1}}
\newcommand{\stX}{\stack{X}}
\newcommand{\stY}{\stack{Y}}

\newcommand{\stwm}{\stack{W}\mspace{-6mu}\stack{M}}

\newcommand{\grstack}[1]{\mathscr{#1}}
\newcommand{\grg}{\grstack{G}}
\newcommand{\grh}{\grstack{H}}
\newcommand{\gre}{\grstack{E}}

\newcommand{\gerbe}[1]{\mathscr{#1}}
\newcommand{\gP}{\gerbe{P}}
\newcommand{\gQ}{\gerbe{Q}}

\newcommand{\bistack}[1]{\mathfrak{#1}}

\newcommand{\CM}{\bicat{XMod}} 
\newcommand{\cm}{\bistack{XMod}} 

\newcommand{\tors}{\operatorname{\textsc{TORS}}}
\newcommand{\bitors}{\operatorname{\textsc{BITORS}}}

\newcommand{\grstacks}{\operatorname{\textsc{Gr\mbox{-}STACKS}}}

\newcommand{\gerbes}{\operatorname{\textsc{GERBES}}}

\title{Butterflies II: Torsors for 2-group stacks}

\author{Ettore Aldrovandi\\
  \small Department of Mathematics,
  Florida State University\\
  \small 1017 Academic Way,
  Tallahassee, FL 32306-4510, USA \\
  \small \url{aldrovandi@math.fsu.edu}%
  \and%
  Behrang Noohi\\
  \small Department of Mathematics,
  King's College London\\
  \small  Strand, London WC2R 2LS, UK\\
  \small \url{behrang.noohi@kcl.ac.uk}}

\date{}

\begin{document}

\maketitle%

\begin{abstract}
  We study torsors over 2-groups and their morphisms.  In
  particular, we study the first non-abelian cohomology group
  with values in a 2-group. Butterfly diagrams encode morphisms
  of 2-groups and we employ them to examine the functorial
  behavior of non-abelian cohomology under change of
  coefficients.  We re-interpret the first non-abelian cohomology
  with coefficients in a 2-group in terms of gerbes bound by a
  crossed module.  Our main result is to provide a geometric
  version of the change of coefficients map by lifting a gerbe
  along the ``fraction'' (weak morphism) determined by a
  butterfly.  As a practical byproduct, we show how butterflies
  can be used to obtain explicit maps at the cocycle level.  In
  addition, we discuss various commutativity conditions on
  cohomology induced by various degrees of commutativity on the
  coefficient 2-groups, as well as specific features pertaining
  to group extensions.
\end{abstract}

\tableofcontents%

\section{Introduction}
\label{sec:introduction}

This paper is the second part of a series aimed at a systematic
study of $n$-group stacks and their torsors.  The first part,
\citep{ButterfliesI}, is dedicated to the case $n=2$ of 2-group
stacks, or gr-stacks, in a slightly older terminology, and
especially their morphisms. The most important result is that if
2-group stacks are made strict by replacing them with (sheaves
of) crossed modules, the groupoid of morphisms between 2-group
stacks is equivalent to that of certain special diagrams called
butterflies between corresponding crossed modules.  This allows one
to overcome the longstanding problem, even present in the non
sheaf-theoretic setting, that replacing a monoidal category with
a strict one is not a functorial construction.

Moving up one step in the cohomological ladder, the present
paper, which is a direct sequel to \citep{ButterfliesI}, is
concerned with the \emph{torsors} for 2-group stacks.  In a very
general sense, torsors are the global geometric objects from
which 1-cocycles with values in a 2-group stack arise, once
suitable local trivializing data have been chosen.  In effect,
after a rigidification has been performed by replacing a 2-group
stack by a crossed module, such cocycles will take values in a
complex of sheaves (of length 2).  This is the categorified
version of the familiar process which associates to a principal
$G$-bundle (or ordinary $G$-torsor) with local sections a
1-cocycle with values in $G$.  Indeed the case $n=1$ is the one
of ordinary group objects.  (In general a similar situation holds
in the case of $n$-group stacks, as we shall see in later
installments of this series.)

Our aim is to the study morphisms of torsors by harnessing the
power of butterflies developed in the first part of this series,
and to illustrate a few applications.

\subsection{Content of the paper}
\label{sec:content-paper}

It is useful to describe the context of our work in general terms.  If
$F\colon \grh\to \grg$ is a morphism of 2-group stacks over a certain
site $\s$, we want an appropriate morphism
\begin{equation}
  \label{eq:1}
  F_*\colon \tors (\grh) \lto \tors (\grg),
\end{equation}
where $\tors (\grg)$ denotes the 2-stack of $\grg$-torsors. One
obtains in this way a geometric definition of degree-one
non-abelian cohomology sets, with built-in functoriality. Namely,
if by $\tors (\grg)(*)$ we denote the 2-groupoid of global
torsors, we can define $\H^1(\grg)$ simply as $\pi_0 \bigl(\tors
(\grg)(*)\bigr)$, the connected components of that 2-groupoid;
once $F_*$ is defined, the functoriality of the first cohomology
follows automatically.

A viable general mechanism by which torsors are extended
``along'' a morphism of $n$-group stacks is in fact well-known:
given an $\grh$-torsor $\stX$, one defines $F_*$ via the
``contracted product''
\begin{equation}
  \label{eq:2}
  F_*(\stX) = \stX\cprod{\grh}\grg,
\end{equation}
see \citep[\S 6]{MR92m:18019}, and
section~\ref{sec:extension-torsors} below for all the details.
The construction on the right-hand side above is the
``categorification'' of the standard one in the case of ordinary
torsors, that is $n=1$.  The above definition of $F_*$ provides a
conceptual answer to finding a morphism~\eqref{eq:1}, and
therefore, by the above geometric definition of cohomology, an
induced morphism
\begin{equation}
  \label{eq:3}
  \H^1(\grh) \lto \H^1(\grg).
\end{equation}
On the other hand, the recently introduced butterfly diagrams
afford a rather fine-grained picture of morphisms of 2-group
stacks, to be recalled below, so one asks for a similar
description of~\eqref{eq:1} and the induced map~\eqref{eq:3}.

To discuss this, let us recall from the first part that a
butterfly allows us to decompose a morphism $F\colon \grh\to \grg$
into a ``fraction''
\begin{equation*}
  \xymatrix@1{%
    \grh & \ar[l]_<<<Q   \gre \ar[r]^P & \grg},
\end{equation*}
where $Q$ is an equivalence of 2-group stacks.  Actually, if we
introduce crossed modules $H_1\to H_0$ and $G_1\to G_0$ for
$\grh$ and $\grg$, respectively, the fraction above is determined
by a \textbf{butterfly diagram} of group objects:
\begin{equation*}
  \xymatrix@R-0.5em{%
    H^{-1}\ar[dd] \ar@/_0.1pc/[dr]  & &
    G^{-1} \ar@/^0.1pc/[dl] \ar[dd]\\
    & E \ar@/_0.1pc/[dl] \ar@/^0.1pc/[dr] &  \\
    H^0 & & G^0
  }
\end{equation*}
The NW-SE sequence is a complex, and the NE-SW sequence is a
group extension.  One finds the resulting map $H_1\times G_1\to
E$ is a crossed module in its own right, which is
quasi-isomorphic to $H_1\to H_0$, and determines the stack
$\gre$.  In sum, with a  butterfly we can split $F\colon \grh\to
\grg$ into a fraction of morphisms corresponding to morphisms of
crossed modules.  In fact the butterfly corresponds to a
fraction in the derived category of crossed modules
\begin{equation}
  \label{eq:4}
  \xymatrix@1{%
    H_\bullet & \ar[l]_<<<q   E_\bullet \ar[r]^p & G_\bullet},
\end{equation}
where now $p$ and $q$ are genuine morphisms of crossed modules
(the latter being a quasi-isomorphism) inducing the corresponding
ones denoted by upper-case letters between corresponding 2-group
stacks.

We have alluded to the fact that classes in, say, $\H^1(\grg)$
can be represented by 1-cocycles with values in the
crossed-module $G_1\to G_0$.  Let us remind the reader, following
\citep{MR92m:18019}, that such cocycles can equivalently be
described as simplicial maps from hypercovers of objects of $\s$
to a reasonable model of the classifying space of $\grg$. One
such is provided, for instance, by the $\W$ construction applied
to the simplicial group $\smp{G}$ determined by the crossed
module. It is possible to prove, using~\eqref{eq:2}, that a
cocycle with values in $H_\bullet$ determines one with values in
$G_\bullet$.  The argument mostly rests on the construction of a
morphism
\begin{equation*}
  \W\smp{H}\lto \W\smp{G}
\end{equation*}
between classifying objects.  (Note, in passing, that this is the
very definition of \emph{weak} morphism of crossed modules in the
set-theoretic case, see \citep{Noohi:weakmaps}.) Unfortunately,
starting from the morphism $F$ as a whole is not very explicit or
constructive, not only because it requires a chosen
rigidification of .the otherwise weak group laws of $\grh$ and
$\grg$, but chiefly because $F$ does not determine a direct
morphism $H_\bullet \to G_\bullet$ between crossed modules.

Our first result is to exploit the butterfly technology to provide a
much more direct approach to computing the morphism~\eqref{eq:3}.  As
explained in section~\ref{sec:push-cohom-class} below, the
morphism~\eqref{eq:3} can be computed by, figuratively speaking,
lifting a 1-cocycle, or equivalently a simplicial map $\eta\colon
U_\bullet\to \W\smp{H}$ along the diagram~\eqref{eq:4}.  More
concretely, one constructs a new simplicial map $\eta'\colon U_\bullet
\to \W\smp{E}$ such that its projection via $q$ is $\eta$ (possibly
after passing to a finer hypercover which will not be notationally
distinguished); in effect $\eta'$ represents the same class as $\eta$,
since $H_\bullet$ and $E_\bullet$ are quasi-isomorphic.  Then the
sought-after morphism is simply obtained by projecting $\eta'$ along
$p$.  Diagrammatically, we have:
\begin{equation*}
  \xymatrix@1{%
    \eta & \ar[l]_<<<q  \eta' \ar[r]^p & \xi}
\end{equation*}
where $\xi$ denotes the resulting simplicial map or 1-cocycle
with values in $G_\bullet$.

This same method, in simpler form, works for 0-cocycles as well, such
as those dealt with in the first part, and it is expected to do so for
higher degree classes in the case the 2-group stacks involved are
symmetric or Picard.

The construction just outlined embodies the general idea that
informs our main result, a novel geometric construction of the
morphism~\eqref{eq:1}.  Starting from the butterfly
decomposition of $F$ we want to decompose $F_*$ as
\begin{equation*}
  \xymatrix@1{%
    \tors (\grh) & \ar[l]_<<<{Q_*} \tors (\gre) \ar[r]^{P_*} &
    \tors (\grg)},
\end{equation*}
where $P_*$ and $Q_*$ are expected to be simpler than $F_*$, since $P$
and $Q$ each arise from a strict morphism.  Moreover, this
decomposition should be such that passing to cohomology classes, or
better yet to representative cocycles, provides a calculation of the
map~\eqref{eq:3} of cohomology sets outlined above.

Now, in practice, we do not implement our program within the context
of torsors over a 2-group stack, essentially due to the fact that the
direction of $P$ is at odds with the natural notion of extension of
torsors along a morphism (\ie $P$ goes in the wrong direction).  One
can of course make the choice of a quasi-inverse $P^*$ to it, but that
defeats the purpose, so to speak; we want something more canonical.

It turns out the concept of gerbes ``bound'' by a crossed module is
the appropriate notion.  In very broad terms, the general idea,
originally due to Debremaeker (see \citep{MR0480515}), is that a gerbe
$\gP$ bound by a crossed module $G_1\to G_0$ is a gerbe equipped with
a morphism
\begin{equation*}
  \mu\colon \gP \lto \tors (G_0)
\end{equation*}
subject to certain additional conditions, recalled in
section~\ref{sec:gerbes-bound-crmod}, which in particular make $\gP$
into a $G_1$-gerbe.  These gerbes give rise to non-abelian cohomology
classes with values in the crossed module (or in fact in a 2-group
stack) too.  Torsors do the same of course, and indeed we prove there
in an equivalence
\begin{equation}
  \label{eq:5}
  \tors (\grg) \lto \gerbes (G_1,G_0),
\end{equation}
which generalizes a similar result of Breen (for the 2-group stack of
$G$-bitorsors for a group object $G$ and $G$-gerbes) put forward in
\citep{MR92m:18019}.  While the equivalence and the statement have
pretty much identical forms, the proof is however quite different, and
we have included it here.

Thus the actual version of the decomposition we provide is is to
\emph{define} a morphism
\begin{equation*}
  F_+ \colon \gerbes(H_1,H_0) \lto \gerbes (G_1,G_0)
\end{equation*}
by means of the following diagram
\begin{equation*}
  \xymatrix@1{%
    \gerbes (H_1,H_0) & \ar[l]_<<<{Q^0_+} \gerbes (E_1,E_0) \ar[r]^{P^0_+} &
    \gerbes (G_0,G_1)}
\end{equation*}
where the definition of $P^0_+$ and $Q^0_+$ is direct (available
in \citep{MR0480515}), since $p$ and $q$ are \emph{strict}
morphisms of crossed modules.  The quasi-inverse to the arrow
pointing to the left, is surprisingly simple in the gerbe
context: from a gerbe $\gQ$ bound by the crossed module
$H_\bullet$, the gerbe bound by $E_\bullet$ that we need is
simply the stack fibered product:
\begin{equation*}
  \gQ'  = \gQ \times_{\tors (H_0)} \tors (E).
\end{equation*}
The image of $\gQ'$ by $Q_+$ is equivalent to $\gQ$, and by
``pushing'' along $P$, that is, considering the image under $P_+$, we
obtain a gerbe bound by $G_\bullet$.  We then prove, essentially by
comparing cohomology classes, that $F_+$ is equivalent to $F_*$,
modulo the equivalence~\eqref{eq:5}, so in other words we obtain a
square
\begin{equation*}
  \xymatrix{%
    \tors (\grh) \ar[d] \ar[r]^{F_*} & \tors (\grg)
    \ar[d] \\
    \gerbes (H_1,H_0) \ar[r]^{F_{+}} & \gerbes (G_1,G_0)
  }
\end{equation*}
commuting up to natural isomorphism.  

After having gone through these general results, we move on to
consider some applications, mainly to the abelian structures on
cohomology resulting when braided, symmetric, or Picard
structures are imposed on the coefficients, and specifically when
group extensions in the sense of Grothendieck
(\citep{MR0354656-VII}) and Breen (\citep[\S 8]{MR92m:18019}) are
concerned.  In the end we make contact with the definition of
weak morphism between crossed modules as simplicial maps between
classifying spaces.  Since several results are already known, our
discussion assumes a more informal character compared to the
previous sections, and many arguments are just sketched.

Let us conclude with a comment about the use of gerbes bound by
crossed modules.  The original intent behind the introduction of the
concept of gerbe bound by a crossed module was to correct the
perceived lack of functoriality inherent in Giraud's definition of
higher non-abelian cohomology using liens (see \citep{MR49:8992}).
Functoriality was addressed in Debremaeker's paper \citep{MR0480515}
by considering only morphisms of crossed modules, that is what we now
call strict morphisms.  This restriction to strict morphisms is not
the natural thing to do, and since non abelian cohomology depends on
the associated 2-group stack, rather than on the coefficient crossed
module itself, introducing torsors led to a better conceptual
understanding of the functoriality of non abelian cohomology.  Thus
the notion has not been developed or used until recently, when it
became useful in different contexts (see for instance
\citep{doi:10.1016/j.jpaa.2007.07.020,math.AG/0301304}).

This state of affairs has been changed by the better control of
morphisms afforded by the use of butterflies, since they allow a
description of all morphisms by way of crossed modules.  Thus now
the use of gerbes bound by crossed modules \emph{plus} the use of
butterflies affords a geometrization of the non abelian derived
category equivalent to the one obtained by using the torsor
picture.

\subsection{Organization of the paper}
\label{sec:organization-paper}

Here is a brief synopsis of this paper's content. Since this is a
direct continuation of \citep{ButterfliesI}, the reader will
unavoidably be constantly referred to that paper. In order to make
this process a little less burdensome, we recall in
section~\ref{sec:recoll-results-from} some of the results of the first
part that we shall most often need here.  In
section~\ref{sec:tors-non-ab-cohom} we have collected results and
definitions concerning torsors over gr-stacks and non-abelian
cohomology.  Our purpose was of course to make a moderate attempt at
being self-contained and at a uniformity of conventions.

By design the material in these sections is not new, except maybe
in the presentation.  New results begin in earnest in
section~\ref{sec:push-cohom-class}, where we explicitly describe
in terms of butterflies the morphism of non-abelian first
cohomology sets induced by a morphism of gr-stacks.

In section~\ref{sec:gerbes-bound-crmod} we present the idea of a
gerbe bound by a crossed module, originally due to Debremaeker.
In addition to re-introduce the main definitions, we analyze the
local structure and prove the cohomology class determined by such
an object takes values in the gr-stack associated to the crossed
module.  Since this is almost the same idea as that of a torsor
for said gr-stack, we determine the precise relation between the
two.  In this way we obtain a generalization of an analogous
result due to \citep[Proposition 7.3]{MR92m:18019}. The sort of
rigidification that the passage from $G\to \Aut (G)$ to a general
crossed module $G\to \Pi$ entails makes the proof very different,
so we discuss it in detail.

The morphism of first non-abelian cohomology sets induced by a
morphism of gr-stacks discussed in purely algebraic terms in
section~\ref{sec:push-cohom-class} has a well-known geometric
realization in terms of extension of torsors along that morphism
(this is the categorification of the well-known extension of
structural groups for principal bundles).  The analogous
procedure in terms of gerbes bound by crossed modules is
described in section~\ref{sec:Extension-butterfly}. It
generalizes Debremaeker's notion of morphism of gerbes bound by
crossed modules, which only uses what we call \emph{strict}
morphisms of crossed modules.  The general case is treated in
section~\ref{sec:extension-butterfly}.  We prove that the
morphism so obtained is equivalent, modulo the equivalence
between torsors and gerbes, to the morphism given by the
extension of torsors, and in section~\ref{sec:map-non-ab} we show
that the induced cohomology class is precisely the one computed
by the procedure described in section~\ref{sec:push-cohom-class}.

Sections~\ref{sec:comm-cond} and~\ref{sec:butt-extens} are
devoted to some applications.  In section~\ref{sec:comm-cond} we
briefly analyze the commutativity conditions on cohomology
ensuing from the assumption that the coefficient crossed module
(or gr-stack) be at least braided.  It is well-known that in this
case the first a priori non-abelian cohomology acquires a group
structure which becomes abelian if the coefficient gr-stack is
symmetric.  Our approach is to analyze these structures in terms
of specific butterfly diagrams associated to braided crossed
modules which express the fact that for a braided gr-stack the
monoidal structure is a weak morphism.  This is discussed in
detail in \citep[\S 7]{ButterfliesI}.  Using these special
butterflies, we are in position to apply the general theory of
section~\ref{sec:Extension-butterfly} to obtain a novel
description of the group structures on cohomology, for which we
can write explicit product formulas at the cocycle level.
Section~\ref{sec:butt-extens} contains some remarks about group
extensions. First about how the classical Schreier theory of
extensions, from the geometric perspective of Grothendieck and
Breen, fits in the butterfly framework.  We then discuss again
commutative structures, and to some extent abelianization maps.
Some final informal paragraphs are devoted to making contact with
the simplicial definition of weak morphism of crossed modules.

\subsection{Conventions and notations}
\label{sec:conv-notat}

In the sequel we shall refer to \citep{ButterfliesI} simply as ``Part
I.''  We keep its standing assumptions, notations, and typographical
conventions: in particular, $\s$ denotes quite generally a site with
subcanonical topology, and $\T=\s\sptilde$ denotes the topos of
$\Set$-valued sheaves over $\s$.  Again as in Part~I we break our
convention usage in the introduction by reverting to the older term
``gr-stack'' in place of the more recent 2-group (stack).  Concerning
the numbering scheme, references to the first part are made using that
paper's numbering sequence.  For this one, we have chosen to cut the
numbering off by one level, due to its reduced length (compared to
\citepalias{ButterfliesI}).

\section{Recollection of results from \citepalias{ButterfliesI}}
\label{sec:recoll-results-from}

\subsection{Crossed modules and gr-stacks}
\label{sec:crossed-modules-gr}

Let $\grg$ be a gr-stack (or 2-group stack), that is a stack over
$\s$ endowed with a group-like monoidal structure
\begin{equation*}
  \otimes \colon \grg\times \grg \lto \grg,
\end{equation*}
see, for example, \cite{MR92m:18019,MR93k:18019,MR95m:18006}, and
\cite{sinh:gr-cats,MR1250465} for the point-wise case. Many of
the results from the previous references which are required in
this text are summarized in \citepalias{ButterfliesI}, to which
the reader is referred for more details. Here we limit ourselves
to recall that starting from $\grg$ we can always construct a
homotopy fibration
\begin{equation*}
  \xymatrix@1{%
    G_1 \ar[r]^\del & G_0 \ar[r]^{\pi_\grg} & \grg
  },
\end{equation*}
where $\del\colon G_1\to G_0$ has the structure of a crossed
module, so that in fact $\grg$ can be recovered as its associated
gr-stack.  More precisely, the crossed module $G_1\to G_0$
provides us with a concrete model for the associated gr-stack,
namely there is an equivalence
\begin{equation*}
  \grg \lisoto \tors(G_1,G_0).
\end{equation*}
Following Deligne \cite{MR546620}, the right-hand side denotes
the stack of those $G_1$-torsors which become trivial after
extension $P\rightsquigarrow P\cprod{G_1}G_0$.  Thus, $\grg$ is
realized as the homotopy fiber
\begin{equation*}
  \xymatrix@1{\grg \ar[r] &
    \tors (G_1) \ar[r]^{\del_*} & \tors (G_0)},
\end{equation*}
where an object of $\grg$ is a pair $(P,s)$, comprising a
right $G_1$-torsor $P$ and a trivialization $s\colon
P\cprod{G_1}G_0\isoto G_0$. When combined with the crossed module
structure, this picture allows us to realize $\grg$ as a
sub-gr-stack of $\bitors (G_1)$ by observing that the underlying
$G_1$-torsor in the pair $(P,s)$ acquires a $G_1$-bitorsor
structure by defining a left $G_1$-action through $s$ as:
\begin{equation*}
  g\cdot p \coloneq p\,g^{s(p)},
\end{equation*}
where $p\in P$, $g\in G_1$, and $s$ is viewed as a
$G_1$-equivariant morphism $s\colon P\to G_0$.  A morphism
$\phi\colon (P,s)\to (Q,t)$ in $\grg$ is therefore a commutative
diagram
\begin{equation*}
  \xymatrix@-.75pc{%
    P \ar[rr]^\phi \ar@/_/[dr]_s && Q \ar@/^/[dl]^t \\ & G_0 &}.
\end{equation*}
It follows that the monoidal structure of $\grg$ can be expressed
through standard contraction of bitorsors: for two objects
$(P,s)$ and $(Q,t)$ of $\grg$ we set
\begin{equation*}
  (P,s)\otimes (Q,t) = (P\cprod{G_1}Q , s\wedge t),
\end{equation*}
where $s\wedge t$ is the $G_1$-equivariant map given by
$(p,q)\mapsto s(p)t(q)$, where $(p,q)$ represents a point of
$P\cprod{G_1}Q$.  It results from the compound trivialization:
\begin{equation*}
  \bigl( P\cprod{G_1} Q\bigr) \cprod{G_1} G_0 \iso 
  P\cprod{G_1} \bigl(  Q \cprod{G_1} G_0 \bigr)
  \xrightarrow{1\wedge t} 
  P\cprod{G_1} G_0 \overset{s}\lto  G_0.
\end{equation*}

In dealing with gr-stacks and crossed modules we will
always---often tacitly---make use of the interplay outlined in
the previous paragraphs, and therefore move freely between
gr-stacks and crossed modules.

\subsection{Butterflies and weak morphisms}
\label{sec:butt-weak-morph}

Let $H_\bullet$ and $G_\bullet$ be crossed modules of $\T$, and
let $\grh$ and $\grg$ denote their associated gr-stacks,
respectively.

A morphism $F\colon \grh\to \grg$, that is, an additive functor,
is by definition a \emph{weak morphism} from $H_\bullet$ to
$G_\bullet$.  All weak morphisms from $H_\bullet$ to $G_\bullet$
form a groupoid, denoted $\cat{WM}(H_\bullet, G_\bullet)$.
  
A \emph{butterfly} from $H_\bullet$ to $G_\bullet$ is by
definition a commutative diagram of group objects of $\T$:
\begin{equation}
  \label{eq:6}
  \vcenter{%
    \xymatrix@R-0.5em{%
      H_1\ar[dd]_\del \ar@/_0.1pc/[dr]^\kappa  & &
      G_1 \ar@/^0.1pc/[dl]_\imath \ar[dd]^\del\\
      & E\ar@/_0.1pc/[dl]_\pi \ar@/^0.1pc/[dr]^\jmath &  \\
      H_0 & & G_0
    }}
\end{equation}
such that the NW-SE sequence is a complex, and the NE-SW sequence
is a group extension.  The various maps satisfy the equivariance
conditions written set-theoretically as:
\begin{equation}
  \label{eq:7}
  \imath (g^{\jmath(e)}) = e^{-1} \imath (g) e,\quad
  \kappa (h^{\pi (e)}) = e^{-1} \kappa (h) e
\end{equation}
where $g\in G_1, h\in H_1, e\in E$.  An easy consequence of
\eqref{eq:7} is that the images of $\jmath$ and $\kappa$ commute
in $E$.

The short-hand notation $[H_\bullet,E,G_\bullet]$ will be used
for a butterfly from $H_\bullet$ to $G_\bullet$.

A \emph{morphism of butterflies} $\phi\colon
[H_\bullet,E,G_\bullet]\to [H_\bullet,E',G_\bullet]$ is given by
a group isomorphism $\phi\colon E\isoto E'$ such that the
diagram:
\begin{equation*}
  \xymatrix@C+1pc{%
    H_1 \ar[r] \ar@/_0.1pc/[dr] \ar[dd] &  E' \ar@/^/[ddr]|(.35)\hole
    \ar@/_/[ddl]|(.35)\hole  & G_1 \ar[l] \ar[dd] \ar@/^0.1pc/[dl] \\
    &  E \ar@/_0.2pc/[dl] \ar@/^0.2pc/[dr] \ar[u] \\
    H_0  & & G_0}
\end{equation*}
commutes and is compatible with all the conditions involved in
diagram~\eqref{eq:6}. Two morphisms are composed in the obvious
way.  In this way butterflies from $H_\bullet$ to $G_\bullet$
form a groupoid, denoted $\cat{B}(H_\bullet, G_\bullet)$.

One of the main results of \citepalias[Theorem
4.3.1]{ButterfliesI} reads, in part:
\begin{theorem}
  \label{thm:1}
  There is an equivalence of groupoids
  \begin{equation*}
    \cat{B}(H_\bullet, G_\bullet) \lisoto \cat{WM}(H_\bullet, G_\bullet).
  \end{equation*}
\end{theorem}
A pair of quasi-inverse functors
\begin{gather*}
  \Phi\colon \cat{B} (H_\bullet, G_\bullet) \lto
  \cat{WM} (H_\bullet, G_\bullet) \\
  \intertext{and} \Psi \colon \cat{WM} (H_\bullet, G_\bullet)
  \lto \cat{B} (H_\bullet, G_\bullet).
\end{gather*}
is explicitly described in Part I.

Strict morphisms of crossed modules (described in detail in Part
I, section 3.2) correspond to butterfly diagrams whose NE-SW
diagonal is split---with a definite choice of the splitting
morphism, see Part I, section 4.5.  Conversely, a
\emph{splittable} butterfly, namely one whose NE-SW diagonal is
in the same isomorphism class as a semi-direct product, by
definition corresponds to a morphism equivalent to a strict one.

A butterfly diagram is called \emph{flippable} or
\emph{reversible} if both diagonal are extensions. The
corresponding weak morphism is an equivalence.

It easy to verify that from the butterfly diagram~\eqref{eq:6}
the homomorphism
\begin{equation*}
  \del_E \colon H_1\times G_1 \lto E,
\end{equation*}
where $\del_E(h,g)= \kappa (h)\imath (g)$, is a crossed module
with the obvious action of $E$ on $H_1\times G_1$ through that of
$H_0$ and $G_0$ on the respective factors. Let us denote this
crossed module by
\begin{equation*}
  E_\bullet\colon E_1\to E_0,
\end{equation*}
with $E_0=E$ and $E_1=H_1\times G_1$.

From Part I we have that the weak morphism given by the
butterfly~\eqref{eq:7} factorizes as a ``fraction''
\begin{equation*}
  \xymatrix@1{%
    H_\bullet & E_\bullet \ar[l]_\sim \ar[r] & G_\bullet}
\end{equation*}
of strict morphisms of crossed modules. The one to the left is a
quasi-isomorphisms, that is, it induces isomorphisms on the
corresponding homotopy sheaves:
\begin{equation*}
  \pi_i (E_\bullet) \iso \pi_i (H_\bullet),\quad i=0,1.
\end{equation*}

\subsection{Composition of butterflies and the bicategory of
  crossed modules}
\label{sec:bicat-cross-modul}

Composition of butterflies is by juxtaposition: Given two
butterflies
\begin{equation*}
  \vcenter{%
    \xymatrix@R-0.5em{%
      K_1\ar[dd]_{\del_K} \ar@/_0.1pc/[dr]  & &
      H_1 \ar@/^0.1pc/[dl]_{\imath'} \ar[dd]^{\del_H}\\
      & F\ar@/_0.1pc/[dl] \ar@/^0.1pc/[dr]^{\jmath'} &  \\
      K_0 & & H_0
    }} \qquad
  \vcenter{%
    \xymatrix@R-0.5em{%
      H_1\ar[dd]_{\del_H} \ar@/_0.1pc/[dr]^\kappa  & &
      G_1 \ar@/^0.1pc/[dl] \ar[dd]^{\del_G}\\
      & E\ar@/_0.1pc/[dl]_\pi \ar@/^0.1pc/[dr] &  \\
      H_0 & & G_0
    }}
\end{equation*}
their composition is the butterfly (defined set-theoretically in
\cite{Noohi:weakmaps}):
\begin{equation*}
  \xymatrix{%
    K_1\ar[dd]_{\del_K} \ar@/_0.1pc/[dr]  & &
    G_1 \ar@/^0.1pc/[dl] \ar[dd]^{\del_G}\\
    & {\displaystyle F\times_{H_0}^{H_1}E}
    \ar@/_0.1pc/[dl] \ar@/^0.1pc/[dr] &  \\
    K_0 & & G_0
  }
\end{equation*}
The center is given by a kind of pull-back/push-out construction:
we take the fiber product $F\times_{H_0}E$ and mod out the image
of $H_1$ (see also \citepalias[\S 5.1]{ButterfliesI}, for
details).

This composition is \emph{not} associative: if $[L_\bullet, M,
K_\bullet]$ is a third butterfly, then the construction of the
composite only yields an \emph{isomorphism}
\begin{equation*}
  \bigl( M\times^{K_1}_{K_0}F \bigr) \times^{H_1}_{H_0} E \lisoto
  M\times^{K_1}_{K_0} \bigl( F \times^{H_1}_{H_0} E \bigr).
\end{equation*}

An almost immediate consequence is
\begin{theorem}[\protect{\citetalias{ButterfliesI}}, Theorem
  5.1.4]
  When equipped with the morphism groupoids $\cat{B}(-,-)$,
  crossed modules in $\cat{T}$ form a bicategory, denoted $\CM
  (\s)$.
\end{theorem}

There are fibered analogs of the various entities we have
introduced so far: so, for instance, one defines a fibered
category $\stack{B} (H_\bullet, G_\bullet)$, which is defined as
usual by assigning to $U\in \Ob\s$ the groupoid
\begin{equation*}
  \cat{B}(H_\bullet\rvert_U,G_\bullet\rvert_U),
\end{equation*}
and to every arrow $V\to U$ of $\s$ the functor
\begin{equation*}
  \cat{B}(H_\bullet\rvert_U,G_\bullet\rvert_U) \lto
  \cat{B}(H_\bullet\rvert_V,G_\bullet\rvert_V).
\end{equation*}
Starting from $\cat{WM}(H_\bullet, G_\bullet)$ instead, an identical
procedure leads to a fibered category $\stwm (H_\bullet, G_\bullet)$
over $\s$.  It is proved in \citepalias[4.6.1, 4.6.2]{ButterfliesI}
that both are stacks (in groupoids) over $\s$. In a more general, but
similar, fashion, the bicategory $\CM (\s)$ has a fibered analog,
denoted $\cm (\s)$. Thanks to the fact that $\stack{B}(H_\bullet,
G_\bullet)$ is itself a stack, $\cm(\s)$ is a pre-bistack over $\s$.
On the other hand, gr-stacks form a 2-stack denoted $\grstacks(\s)$,
hence the obvious morphism $\cm (\s) \to \grstacks (\s)$ sending a
crossed module to its associated gr-stack is 2-faithful.  Moreover,
every gr-stack $\grg$ is equivalent to the gr-stack associated to a
crossed module---see \citepalias[Proposition
5.3.7]{ButterfliesI}. Therefore the above morphism is essentially
surjective, and it follows that $\cm (\s)$ is a bistack.

\section{Torsors and non-abelian cohomology}
\label{sec:tors-non-ab-cohom}

In this section we recall some facts about $\grg$-torsors, where
$\grg$ is a gr-stack.  This is necessary in order to compare them with
one of the main objects of study in this text, the gerbes bound by the
crossed module $G_1\to G_0$ whose associated gr-stack is $\grg$.
Those gerbes will be introduced in
section~\ref{sec:gerbes-bound-crmod}.  Since we shall also be
concerned with classes of equivalence of such objects, as well as
functoriality properties, it is useful to go through a quick review of
some definitions in non-abelian cohomology.

\subsection{Non-abelian cohomology}
\label{sec:non-abel-cohom}

Let us recall the main definitions, following \cite{MR92m:18019}
and \cite{MR0491680,MR991977,MR862637}.  Let $\smp{G}$ be a
simplicial group-object of $\T$.  The non-abelian cohomology with
values in $\smp{G}$ can be defined as
\begin{equation*}
  \H^i (*,\smp{G}) =
  \begin{cases}
    \Hom_{\mathcal{D}(\T)} (*,
    \operatorname{\Omega}^{-i}\smp{G}),  &  i\leq 0,\\
    \Hom_{\mathcal{D}(\T)} (*, \B\smp{G}), & i=1.
  \end{cases}
\end{equation*}
Here $*$ denotes the terminal object of $\T$, $\Omega$ denotes
the loop construction, whereas $\B{\smp{G}}$ is some (in fact
any) form for the classifying space construction, for example
$\W\smp{G}$.  $\mathcal{D}(\T)$ denotes the derived category of
simplicial objects of $\T$ in the same sense as
\cite{MR0491680,MR92m:18019}, that is, by localizing at the
morphisms of simplicial objects that induce isomorphisms of
homotopy sheaves.

Note that the simplicial group structure is only relevant in
order to define $\H^1$, whereas for all other degrees $i\leq 0$
the definition only uses the underlying simplicial set
structure. But also note that the former will only be a pointed
set, as opposed to the others which carry group structures
(abelian for $i<0$).  If we use the convention that
$\B^{-1}\eqdef \Omega$, the various $\H^i(*,\smp{G})$ are
computed as a colimit:
\begin{equation*}
  \H^i(*,\smp{G}) = \colim_{V\to *}\bigl[ *, \B^i\smp{G} \bigr],
\end{equation*}
where the colimit runs over homotopy classes of hypercovers of
$*$ and $[-,-]$ denotes (simplicial) homotopy classes.

Our main focus will be the pointed set $\H^1(*,\smp{G})$ when the
coefficient simplicial group arises from a crossed module $G_1\to
G_0$, which we denote by $\H^1(*, G_1\to G_0)$.  In view of the
fact that any gr-stack $\grg$ can be realized as the gr-stack
associated to a crossed module $G_1\to G_0$, as explained in
Part~I, we can write the same cohomologies by emphasizing the
stack, rather than the crossed module, as coefficients, as
$\H^i(*,\grg)$, $i\leq 1$.  In fact more stress will be put on
the \emph{cocycles} representing cohomology classes, rather than
on the classes themselves. After all, the former naturally arise
from any appropriate decomposition (\ie local description) of
geometric objects, such as torsors and gerbes, as it will be
clear below.

Following \cite{MR92m:18019}, it will be convenient to recall the
simplicial definition of 1-cocycles, as well as the more
geometric one that simply categorifies the standard definition by
replacing a group with a gr-stack.

\subsection{1-Cocycles with values in crossed modules}
\label{sec:cocycles}

If $\smp{G}_\bullet$ is a simplicial group object of $\T$, there
is a model for its classifying space provided by the
$\W$-construction.  Namely, $\W\smp{G}_\bullet$ is the simplicial
object of $\T$ given by:
\begin{equation*}
  \W \smp{G}_0 = *, \qquad
  \W \smp{G}_n = \smp{G}_0\times \smp{G}_1\times \dotsb
  \times \smp{G}_{n-1}\,,\quad n\geq 1.
\end{equation*}
The face and degeneracy maps are:
\begin{align*}
  d_i ( \smp{g}_0,\dotsc, \smp{g}_{n-1}) & =
  \begin{cases}
    (d_1\smp{g}_1,\dotsc, d_{n-1}\smp{g}_{n-1}) & i=0 \\
    (\smp{g}_0,\dotsc, \smp{g}_{i-1}d_0\smp{g}_i, \smp{g}_{i+1},
    \dotsc, d_{n-i-1}\smp{g}_{n-1}) & 0<i<n \\
    (\smp{g}_0,\dotsc, \smp{g}_{n-2}) & i=n
  \end{cases} \\
  \intertext{and} s_i(\smp{g}_0,\dotsc, \smp{g}_{n-1} ) & =
  \begin{cases}
    (\smp{1}, s_0\smp{g}_0,\dotsc, s_{n-1}\smp{g}_{n-1} ) & i=0\\
    (\smp{g}_0,\dotsc, \smp{g}_{i-1}, \smp{1},
    s_0\smp{g}_i,\dotsc, s_{n-i-1}\smp{g}_{n-1}) & 0<i<n \\
    (\smp{g}_0,\dotsc, \smp{g}_{n-1},\smp{1}) & i=n
  \end{cases}
\end{align*}
We have slightly changed the formulas of ref.\
\cite[\S{21}]{MR1206474} in order to better fit with our ``action
on the right'' convention.

If $G$ is a group object of $\T$, identified with the constant
simplicial group, then the previous construction reduces to the
standard classifying simplicial space $\B G$.
\begin{definition}
  \label{def:1}
  Let $V_\bullet\to U$ be a hypercover. A \emph{$1$-cocycle over
    $U$} is a simplicial map $\xi\colon V_\bullet\to
  \W\smp{G}_\bullet$.  Two such cocycles $\xi,\xi'$ are
  \emph{equivalent} if there is a simplicial homotopy $\alpha
  \colon \xi \Rightarrow \xi'\colon V_\bullet \to
  \W\smp{G}_\bullet$.
\end{definition}

Let $\smp{G}_\bullet$ be the nerve of the groupoid $\grpd{G}$
determined by a crossed module $G_1\to G_0$.  In this case we
have $\W\smp{G}_1 = G_0, \W\smp{G}_2 = G_0 \times (G_0\times
G_1), \W\smp{G}_3 = G_0 \times (G_0\times G_1)\times (G_0\times
G_1\times G_1)$, etc.  A simplicial map $\xi\colon V_\bullet\to
\W\smp{G}_\bullet$ will be determined by its $3$-truncation
(\cite{MR92m:18019}).

A rather tedious, but otherwise straightforward calculation shows
that the simplicial map $\xi$ determines, and is determined by, a
pair $(x,g)$ where $x\colon V_1 \to G_0$ and $g\colon V_2\to G_1$
satisfying the condition
\begin{subequations}
  \label{eq:8}
  \begin{align}
    \label{eq:9}
    d_1^* x & = d_2^*x\, d_0^*x\, \del g \\
    \label{eq:10}
    d_0^*g\, d_2^*g & = (d_3^*g)^{(d_0d_1)^*x} \, d_1^*g
  \end{align}
\end{subequations}
and the normalizations $s_0^*x=1$, $s_0^*g = s_1^*g = 1$.  The
explicit expressions of the maps $\xi_i$, $i=0,\dotsc,3$ are as
follows: $\xi_0 = *, \xi_1=x\colon V_1\to G_0$, whereas $\xi_2
\colon V_2 \to G_0 \times (G_0 \times G_1)$ and $\xi_3\colon V_3
\to G_0 \times (G_0\times G_1)\times (G_0\times G_1\times G_1)$
are given by
\begin{align*}
  \xi_2 & = ( d_2^*x, (d_0^*x, g)) \\
  \xi_3 & = ( (d_2d_3)^*x, ( (d_0d_3)^*x, d_3^*g), ( (d_0d_1)^*x,
  d_0^*g, (d_0^*g)^{-1} d_1^*g).
\end{align*}
\begin{remark}
  \label{rem:1}
  There exists a compelling way of organizing the above data.
  The idea is that from the form of $\xi_1$ and $\xi_2$ we can
  use $x$ as a label for a 1-cell of $\W\smp{G}_\bullet$, and $g$
  as a label for the 2-cells.  With this in mind,
  equation~\eqref{eq:9} represents a 2-cell with its boundary, as
  in the following diagram:
  \begin{equation*}
    \xymatrix{%
      0 \ar[rr]^{d^*_{2}x} \ar@/_/[dr]_{d^*_{1}x}
      \ar@{}[drr]|g & &
      1 \ar@/^/[dl]^{d^*_{0}x} \\ & 2 &
    }
  \end{equation*}
  Similarly, \eqref{eq:10} represents the compatibility of the
  four possible pullbacks of~\eqref{eq:9}, and therefore has a
  tetrahedral shape:
  
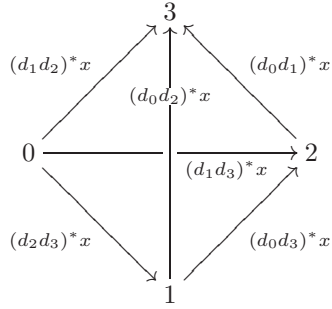
\begin{equation*}
    \xymatrix@ur@+3pc{%
      0 \ar[r]^{(d_1d_2)^*x} \ar[d]_{(d_2d_3)^*x} \ar[dr]|\hole_(0.7){(d_1d_3)^*x} & 3 \\
      1 \ar[r]_{(d_0d_3)^*x} \ar[ru]|(0.7){(d_0d_2)^*x} & 2 \ar[u]_{(d_0d_1)^*x}
    }
  \end{equation*}
  We have not recorded the face labels to avoid cluttering the
  diagram. To recover them, and hence equation~\eqref{eq:10},
  observe that for $i\in \mathbf{3}=\lbrace 0,\dots,3\rbrace$,
  $d_i^*g$ is the 2-cell with vertices given by the complement of
  $i$ in $\mathbf{3}$.

  Alternatively, the following planar version is perhaps clearer:
  \begin{equation*}
    \vcenter
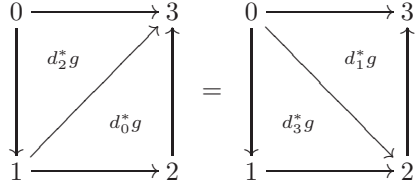
{%
      \xymatrix@+2pc{%
        0 \ar[r] \ar[d] \ar@{}[dr]|(0.3){d_2^*g}|(0.7){d_0^*g} & 3 \\
        1 \ar[r] \ar[ru] & 2 \ar[u]
      } 
    } \; =\;
    \vcenter{%
      \xymatrix@+2pc{%
        0 \ar[r] \ar[d] \ar[dr] & 3 \\
        1 \ar[r] \ar@{}[ur]|(0.3){d_3^*g}|(0.7){d_1^*g}  & 2 \ar[u]
      }
    }
  \end{equation*}
  Note also that the 2-cell $d_3^*g$ is the only one \emph{not}
  including the vertex $3$. Hence an action by $(d_0d_1)^*x$ is
  required.  Also, the right action should match composition, so
  that the 2-cells should be traversed from bottom to top,
  relative to the last diagram.

  Some aspects of the above constructions, in particular the
  seemingly arbitrary labeling of the vertices, given that
  $\W\smp{G}_\bullet$ has only one 0-cell, may appear somewhat
  arbitrary.  A full geometric explanation will be possible after
  the connection with trivializations of torsors and
  (equivalently) gerbes is made in
  sections~\ref{sec:torsors-gr-stacks}
  and~\ref{sec:cohom-class-gerbe}, respectively.
\end{remark}

A \textbf{simplicial homotopy} $\alpha\colon \xi\to \xi'$ is
uniquely determined by $y\colon V_0\to G_0$ and $a_0,a_1\colon
V_1\to G_1$ such that:
\begin{equation}
  \label{eq:11}
  \begin{aligned}
    (d_1^*y)\, x' & = x\, (d_0^*y)\, \del (a_1a_0^{-1}) \\
    d_0^* (a_1a_0^{-1})\,d_2^* (a_1a_0^{-1})^{d_0^*x'} g' &=
    g^{(d_0d_1)^*y}\, d_1^*(a_1a_0^{-1})
  \end{aligned}
\end{equation}
Note that the change $a_0\to a_0a, a_1\to a_1a$ gives another
homotopy between $\xi$ and $\xi'$.

The simplicial homotopy itself (again as in \cite[\S
5]{MR1206474}) in this case is given by maps $\alpha^0_0\colon
V_0\to G_0$, $\alpha^1_i\colon V_1\to G_0\times (G_0\times G_1)$
for $i=0,1$, and $\alpha^2_i\colon V_2\to G_0\times (G_0\times
G_1)\times (G_0\times G_1\times G_1)$, $i=0,1,2$, given by
\begin{align*}
  \alpha^0_0 & = y \\
  \alpha^1_0 & = (d_1^*y, (x', a_0)) \\
  \alpha^1_1 & = (x, (d_0^*y,  a_1)) \\
  \alpha^2_0 & = ((d_1d_2)^*y, (d_2^*x',d_2^*a_1), (d_0^*x',g',
  {g'}^{-1}(d_2^*a_0^{-1})^{d_0^*x}g'\,d_1^*a_0))\\
  \alpha^2_1 & = (d_2^*x, ((d_0d_2)^*y,),(d_0^*x', d_0^*a_0,
  d_0^*a_0^{-1} (d_2^* a_0^{-1})^{d_0^*x'} g' d_1^*a_0  ))\\
  \alpha^2_2 & = (d_2^*x,(d_0^*x, g),((d_0d_1)^*y, d_0^*a_1,
  d_0^*a_1^{-1}d_1^*a_1)\\
\end{align*}
These results are essentially the same (barring a different set
of conventions) as those of \cite[\S 6.4--6.5]{MR92m:18019} for
the crossed module $\iota \colon G\to \Aut (G)$.

\subsection{Bitorsor cocycles}
\label{sec:bitorsor-cocycles}

Let $\grg$ be a gr-stack.  Let $U_\bullet$ be a hypercover, for
example the \Cech complex $\cech U$ of a generalized cover $U\to
*$.
\begin{definition}
  \label{def:2}
  A \emph{1-cocycle} with values in $\grg$ consists of a pair
  $(g,\gamma)$, where $g$ is an object of $\grg$ over $U_1$, and
  $\gamma$ a morphism of $\grg$ over $U_2$, satisfying the
  cocycle conditions
  \begin{subequations}
    \begin{gather}
      \label{eq:12}
      \gamma \colon d_1^*g \lisoto d_2^*g\cdot d_0^*g \\
      \intertext{over $U_2$, and the coherence condition}
      \label{eq:13}
      \bigl((d_2d_3)^*g\cdot d_0^*\gamma \bigr) \circ d_2^*\gamma
      = a\circ \bigl(d_3^*\gamma \cdot (d_0d_1)^*g \bigr) \circ
      d_1^*\gamma,
    \end{gather}
  \end{subequations}
  over $U_3$, where $a$ is the associator isomorphism for the
  group law in $\grg$.  Two cocycles $(g,\gamma)$ and
  $(g',\gamma')$ (assumed for simplicity to be defined over the
  same $U_\bullet$) are \emph{equivalent} if there is a pair
  $(h,\eta)$, where $h\in \Ob\grg_{U_0}$ and $\eta \in \Mor
  \grg_{U_1}$, such that:
  \begin{subequations}
    \label{eq:14}
    \begin{equation}
      \eta \colon g\cdot (d_0^*h) \lisoto (d_1^*h)\cdot g' \\
    \end{equation}
    and the diagram
    \begin{equation}
      \vcenter{%
        \xymatrix@d{%
          (d_1d_2)^*h\cdot d_1^*g' \ar[r]^{\gamma'} &
          (d_1d_2)^*h \cdot (d_2^*g'\cdot d_0^*g') &
          ((d_1d_2)^*h \cdot d_2^*g')\cdot d_0^*g' \ar[l]_a &
          (d_2^*g\cdot (d_0d_2)^*h )\cdot d_0^*g' \ar[d]^a
          \ar[l]_{d_2^* \eta}
          \\
          d_1^*g\cdot (d_0d_1)^*h \ar[u]^{d_1^* \eta} \ar[r]^{\gamma} &
          (d_2^*g\cdot d_0^*g)\cdot (d_0d_1)^*h \ar[r]^a &
          d_2^*g\cdot (d_0^*g\cdot (d_0d_1)^*h) &
          d_2^*g\cdot ((d_0d_2)^*h\cdot d_0^*g') \ar[l]_{d_0^* \eta}
        }
      }
    \end{equation}
    commutes.
  \end{subequations}
\end{definition}
In view of the discussion on the relationship between $\grg$ and
the crossed module reviewed in
sect.~\ref{sec:crossed-modules-gr}, whereby the monoidal
structure of $\grg$ is described in terms of contracted products
of $G_1$-bitorsors, a 1-cocycle such as $(g,\gamma)$ in
Definition~\ref{def:2} will be referred to, albeit imprecisely,
as bitorsor cocycle.

It is easy to pass from a 1-cocycle with values in $\grg$ to a
1-cocycle with values in $\W\smp{G}_\bullet$. Indeed, recall from
\cite{MR92m:18019} or from the remarks in
sect.~\ref{sec:crossed-modules-gr} that $\grg \iso \tors
(G_1,G_0)$, the gr-stack of $G_1$-torsors equipped with a chosen
trivialization of their extensions to $G_0$. Thus $g\in \Ob
\grg_{U_1}$ can be thought of as such an object. In other words,
we may write $g$ as the pair $g=(E,s)$, where $E$ is the
underlying $G_1$-torsor and $s\colon E\to G_0$ is the equivariant
morphism providing the trivialization as a $G_0$-torsor.  So we
have:
\begin{lemma}
  \label{lem:1}
  There exists a refinement $V_\bullet$ of $U_\bullet$ such that
  the bitorsor cocycle $(g,\gamma)$ determines a 1-cocycle
  $V_\bullet\to \W\smp{G}_\bullet$.
\end{lemma}
\begin{proof}
  Let $V\to U_1$ be a generalized cover such that the restriction
  of the underlying $G_1$-torsor $E$ of $g=(E,s)$ becomes
  trivial.  Then by \cite[V, Théorème 7.3.2]{MR0354653} there
  exists a hypercover $V_\bullet$ and a map $V_\bullet \to
  U_\bullet$ which for degree $n=1$ factorizes through the chosen
  cover:
  \begin{equation*}
    V_1 \to V\to U_1.
  \end{equation*}
  Over $V_1$ we have $E\rvert_{V_1}\iso G_1\rvert_{V_1}$, and $s$
  is determined by its value $s(1)\in G_0$.  Thus $g$ may simply
  be identified with this element of $G_0(V_1)$.  In turn, the
  morphism $\gamma$ is identified with an element of $G_1$ over
  $V_2$, since the underlying map of $G_1$-torsors is a morphism
  of trivial torsors.  That is, the required element is simply
  $\gamma (1)\in G_1 (V_2)$.  Notice that from the identification
  of $g$ with $s(1)\in G_0$ it follows that $d_2^*g\cdot d_0^*g$
  is identified with the product $d_2^*s(1) d_0^*s(1)$.  Since
  $\gamma$ is a morphism of $(G_1,G_0)$-torsors, we must have
  that
  \begin{equation*}
    d_1^* s (1) = d_2^*s(\gamma(1)) d_0^*s(\gamma(1))
    = d_2^*s(1) d_0^*s(1) \, \del \gamma (1).
  \end{equation*}
  Furthermore, it is not difficult to realize that the coherence
  condition for $\gamma$ on $V_3$ becomes
  \begin{equation*}
    d_0^*\gamma (1) d_2^*\gamma (1) =
    d_3^*\gamma (1)^{(d_0d_1)^* s(1)}
    d_1^*\gamma (1).
  \end{equation*}
  These are precisely the cocycle relations~\eqref{eq:8} (modulo
  exchanging $x \leftrightarrow g$ and $g\leftrightarrow \gamma$
  in the notation).
\end{proof}
The procedure in the proof of Lemma~\ref{lem:1} will repeatedly
be used in the sequel.
\begin{remark}
  \label{rem:2}
  If $U_\bullet \to \W\smp{G}_\bullet$ is a simplicial map, where
  again $\smp{G}_\bullet$ is the simplicial group determined by a
  crossed module, a converse procedure allows one to obtain a
  1-cocycle with values in the associated gr-stack $\grg$ relative to
  the \Cech nerve $\cech (U_0)$, where $U_0$ is the degree $n=0$
  object of $U_\bullet$.  The (long) proof can be extracted from
  \cite[\S 6.5]{MR92m:18019}. No explicit use will be made of such
  procedure in the rest of this paper.
\end{remark}

\subsection{Torsors for gr-stacks}
\label{sec:torsors-gr-stacks}

The definition of \emph{torsor} under a gr-stack has been given
in full generality in~\cite[6.1]{MR92m:18019}, so here we will
confine ourselves to only recalling the main points.  Let $\grg$
be a gr-stack over $\s$.  In modern parlance, a $\grg$-torsor is
the categorification of the standard notion of torsor, as
follows.

A right-action of $\grg$ on a stack in groupoids $\stX$ is given
by a morphism of stacks
\begin{equation*}
  m \colon \stX \times \grg \lto \stX
\end{equation*}
plus a natural transformation
\begin{equation}
  \label{eq:15}
  \vcenter{%
    \xymatrix@C+1pc{%
      \stX\times \grg\times \grg \ar[r]^{(m,\id_\grg)}
      \ar[d]_{(\id_\stX,\otimes_\grg)} &
      \stX \times \grg \ar[d]^m_{}="m1" \\
      \stX\times \grg \ar[r]_m^{}="m2"
      \ar@/_/@{=>} "m1";"m2" _{\mu}
      & \stX
    }}
\end{equation}
which amounts, for objects $x,g_0, g_1$, to a functorial
isomorphism
\begin{equation*}
  \mu_{x,g_0,g_1}\colon (x\cdot g_0)\cdot g_1
  \lisoto x\cdot (g_0\cdot g_1),
\end{equation*}
where $x\cdot g$ stands for $m(x,g)$. We require that:
\begin{enumerate}
\item \label{item:1} the pair $(m,\mu)$ satisfy the standard
  pentagon diagram;
\item \label{item:2} the composite
  \begin{equation*}
    \xymatrix@1{%
      \stX \ar[r]^(.4)\sim &
      \stX \times \mathbf{1} \ar[r] &
      \stX\times \grg \ar[r]^(.6)m & \stX
    }
  \end{equation*}
  be isomorphic to the identity functor of $\id_\stX$, where
  $\mathbf{1}\to \grg$ sends the unique object to the identity object
  of $\grg$. Moreover, this morphism must be compatible with $m$ and
  $\mu$, in the sense that the two diagrams
  \cite[(6.1.4)]{MR92m:18019}, resulting from combining it
  with~\eqref{eq:15}, must be commutative.
\end{enumerate}
Most importantly, we require that the morphism
\begin{equation*}
  \tilde m = (\mathrm{pr}_1,m) \colon
  \stX \times \grg \lto \stX \times \stX
\end{equation*}
be an equivalence.  Having so far defined what ought to be called
a \emph{pseudo-}torsor, we need to complete the definition by
adding the condition that there exist a (generalized) cover $U\to
*$ such that the fiber category $\stX_U$ be non-empty.

A \emph{morphism} of $\grg$-torsors $\stX \to \stX'$ consists of
a stack morphism $F\colon \stX\to \stX'$ together with a natural
transformation
\begin{equation*}
  \xymatrix@C+1pc{%
    \stX \times \grg \ar[d]_m \ar[r]^{(F,\Id_\grg)} &
    \stX \times \grg \ar[d]^{m'}_{}="m'" \\
    \stX \ar[r]_F^{}="F" & \stX'
    \ar @/_0.5pc/ @{=>} "m'";"F" _\phi
  }
\end{equation*}
compatible with the transformations $\mu$ and $\mu'$ (That is,
with the diagrams~\eqref{eq:15}).

A \emph{2-morphism} of $\grg$-torsors is a 2-morphism $\alpha
\colon F \Rightarrow F'$ such that the diagrams
\begin{equation*}
  \xymatrix@C+1pc{%
    \stX \times \grg \ar[d]_m \ar[r]^{(F,\Id_\grg)} &
    \stX \times \grg \ar[d]^{m'}_{}="m'" \\
    \stX \ar[r]_(0.4)F^{}="F"
    \rlowertwocell_{F'}{\alpha}
    \ar @/_0.5pc/ @{=>} "m'";"F" _\phi & \stX'
  }
  \qquad
  \xymatrix@C+1pc{%
    \stX \times \grg \ar[d]_m \ar[r]_{(F',\Id_\grg)}
    \ruppertwocell^{(F,\Id_\grg)}{\alpha} &
    \stX \times \grg \ar[d]^{m'}_{}="m'" \\
    \stX \ar[r]_{F'}^{}="F'" & \stX'
    \ar @/_0.5pc/ @{=>} "m'";"F'" _{\phi'}
  }
\end{equation*}
define a commutative diagram of 2-morphisms.
\begin{remark}
  We have defined the notion of right torsors. That of left
  torsor is defined in the same way. It is actually the one
  adopted in \cite{MR92m:18019}.
\end{remark}

With the notions of morphism and 2-morphism outlined above,
$\grg$-torsors comprise a 2-category. In fact, all together they
form a neutral 2-gerbe over $\s$ denoted $\tors (\grg)$. The
fiber above $U\in \Ob (\s)$ is the 2-category of
$\grg\rvert_U$-torsors (cf.\ \cite{MR93k:18019,MR95m:18006}).

\subsection{Contracted product of gr-stacks}
\label{sec:prgr_10}

We will need to consider the notion of contracted product of
torsors over a gr-stack in some detail. It is introduced
in~\cite[\S 6.7]{MR92m:18019} (credited to J.~Bénabou). (We use a
slightly different convention for some of the diagrams.)
  
If $\stX$ (resp.\ $\stY$) is a right (resp.\ left) $\grg$-torsor,
or more generally a stack with a $\grg$-action, the
\emph{contracted product} $\stX\cprod{\grg}\!  \stY$ is defined
as follows. The objects are pairs $(x,y)\in \Ob \stX\times \stY$.
A morphism $(x,y)\to (x',y')$ is an equivalence classes of
triples $(a,g,b)$, where $g\in \Ob \grg$, and $a\colon x \to
x'\cdot g$ and $b\colon g\cdot y\to y'$ are morphisms of $\stX$
and $\stY$, respectively. Two triples $(a,g,b)$ and $(a',g',b')$
are equivalent if there is a morphism $\gamma\colon g\to g'$ in
$\grg$ such that the diagrams
\begin{equation*}
  \xymatrix@R-1pc{%
    & x'\cdot g \ar[dd]^{\id_{x'}\cdot\gamma} \\
    x \ar@/^0.3pc/ [ur]^a \ar@/_0.3pc/[dr]_{a'} \\
    & x'\cdot g'
  }\qquad
  \xymatrix@R-1pc{%
    g\cdot y \ar@/^0.3pc/[dr]^b \ar[dd]_{\gamma\cdot \id_y} \\
    & y' \\
    g'\cdot y \ar@/_0.3pc/[ur]_{b'}
  }
\end{equation*}
commute. The composition of two morphisms $(x_1,y_1)\to
(x_2,y_2)$ and $(x_2,y_2)\to (x_3,y_3)$ represented by triples
$(a,g,b)$ and $(a',g',b')$, respectively, is represented by the
triple given by the expected compositions
\begin{gather*}
  x_1 \overset{a}{\lto} x_2\cdot g \overset{a'\cdot g'}{\lto}
  (x_3\cdot g')\cdot g \lisoto
  x_3 \cdot (g'\cdot g) \\
  (g'\cdot g) \cdot y_1 \lisoto g' \cdot (g\cdot y_1)
  \overset{g'\cdot b}{\lto} g'\cdot y_2 \overset{b'}{\lto} y_3
\end{gather*}
and, of course, $g'\cdot g$.

It should be observed that the foregoing procedure produces a
fibered category over $\s$ with group law. We denote by
$\stX\cprod{\grg}\! \stY$ the associated stack.  One may also
characterize $\stX\cprod{\grg}\! \stY$ as the ``2-Limit'' of the
diagram
\begin{equation*}
  \xymatrix@1{%
    \stX \times \stY \times \grg
    \ar@<0.5ex>[r] \ar@<-0.4ex>[r] & \stX \times \stY
  }
\end{equation*}
where one arrow is the projection and the other is the (right)
action $(x,y,g) \to (x\cdot g, g^*\cdot y)$, where $x,y,g$ are
objects and $g^*$ is a choice for the inverse of $g$.

Properties analogous to the familiar ones for ordinary torsors
hold. For example, whereas in the ordinary contracted product
$P\cprod{G}Q$ of $G$-spaces one has the relation
\begin{equation*}
  (xg,y) = (x,gy),
\end{equation*}
namely the two pairs $(xg,y)$ and $(x,gy)$ represent the same point of
$P\cprod{G}Q$, here one has the isomorphism
\begin{equation*}
  (x\cdot g, y) \lisoto (x,g\cdot y)\,,
\end{equation*}
represented by the triple $(\id_{x\cdot g},g,\id_{g\cdot y})$.

\subsection{Cohomology classes and classification of torsors}
\label{sec:cohom-classification}

\begin{proposition}[\protect{\cite[Proposition
    6.2]{MR92m:18019}}]
  \label{prop:1}
  Let $G_1\to G_0$ be a crossed module of $\T$. The elements of
  the pointed set $\H^1(*, G_1\to G_0)$ are in bijective
  correspondence with equivalence classes of right $\grg$-torsors
  over $\s$, where $\grg = \bigl[ G_1\to G_0\bigr]\sptilde$.
\end{proposition}
\begin{proof}[General idea of the proof]
  The central argument goes through the standard computation with
  1-cocycles subordinated to hypercovers $U_\bullet$.  Suppose
  $\stX$ is a right $\grg$-torsor over $\s$, as described above.
  The choice of an object $x$ of $\stX$ over $U_0$ leads to
  establishing the existence of an object $g$ of $\grg$ over
  $U_1$ such that
  \begin{equation*}
    d_0^*x \lisoto d_1^*x \cdot g.
  \end{equation*}
  After pulling back to $U_2$, from the local equivalence of $\stX$
  and $\grg$ we can conclude that there must exist a
  morphism~\eqref{eq:12} over $U_2$, with $\gamma$
  satisfying~\eqref{eq:13} over $U_3$.  The choice of another
  object $x'$ of $\stX$, still over $U_0$ say, leads to another
  1-cocycle $(g',\gamma')$ \emph{equivalent} to $(g,\gamma)$, in
  the sense of Definition~\ref{def:2}; that is, there is a pair
  $(h,\eta)$ where $h$ is an object of $\grg$ over $U_0$ and a
  $\eta$ morphism over $U_1$ satisfying equations~\eqref{eq:14}.

  From a 1-cocycle $(g,\gamma)$ one can extract a 1-cocycle with
  values in the crossed module $G_1\to G_0$ as explained at the
  end of sect.~\ref{sec:bitorsor-cocycles}.

  Conversely, as mentioned in Remark~\ref{rem:2}, the procedure
  from the proof of \cite[Proposition 6.2]{MR92m:18019}, in
  particular \S 6.5, allows us to reconstruct a bitorsor cocycle,
  and ultimately a $\grg$-torsor, from a 1-cocycle with values in
  $G_1\to G_0$.
\end{proof}

\section{Pushing cohomology classes along butterflies}
\label{sec:push-cohom-class}

Changing the coefficients results in a morphism in non-abelian
cohomology.  From the point of view of the general definition
recalled in sect.~\ref{sec:non-abel-cohom}, this is done by means
of a morphism of simplicial groups $\smp{H}_\bullet \to
\smp{G}_\bullet$, which in our case is the one induced by a
morphism of crossed modules, and ultimately by a morphism
$F\colon \grh \to \grg$ of gr-stacks.  We are also specifically
interested in the case $i=1$, and we want to provide a short
account of how the morphism
\begin{equation*}
  F_*\colon \H^1(*,\grh) \lto \H^1(*,\grg)
\end{equation*}
can be profitably described in terms of butterflies. This is a
necessary stepping stone in the more geometric description of the
first non-abelian cohomology group with values in a gr-stack to
be presented further down in the paper.  After some general
observations, we begin with an elementary approach to the above
morphism in terms of explicit 1-cocycles with values in crossed
modules.  We then show how the more conceptual formulation in
terms of bitorsor cocycles can be reduced to these explicit
calculations.

\subsection{General remarks}
\label{sec:general-remarks}

Let $(F,\lambda)\colon \grh \to \grg$ be a morphism of gr-stacks
over $\s$, where we have explicitly marked the natural
isomorphism $\lambda$ providing the additivity:
\begin{equation*}
  \lambda_{y_1,y_2}\colon F (y_1y_2) \lisoto F(y_1)F(y_2),
\end{equation*}
for any two objects $y_1,y_2$ of $\grh$.  The following is an
easy claim whose proof is left to the reader.
\begin{lemma}
  \label{lem:2}
  Let $(F,\lambda)$ be as above, and let $(y,h)$ be a 1-cocycle
  with values in $\grh$ relative to a hypercover $U_\bullet\to *$
  as in Definition~\ref{def:2}.  Then $(F(y),\lambda\circ F(h))$
  is a 1-cocycle with values in $\grg$ (relative to the same
  hypercover).  If $(y,h)$ and $(y',h')$ are two equivalent
  1-cocycles with values in $\grh$, then so are their images
  $(F(y),\lambda\circ F(h))$ and $(F(y'),\lambda\circ F(h'))$.
\end{lemma}
Our goal is to explicitly calculate $(F(y), \lambda\circ F(h))$
by means of a butterfly representing $F$.

\subsection{Lift of a 1-cocycle along a butterfly}
\label{sec:lift-cocycle-along}

Since a butterfly $[H_\bullet, E, G_\bullet]$ corresponds to a
morphism $F\colon \grh \to \grg$, it is expected that it will be
possible to ``lift'' a 1-cocycle $\eta = (y,h)$ with values in
$\W \smp{H}_\bullet$ to one with values in $\W \smp{G}_\bullet$.
Note that, after having observed that the butterfly $E$ or
equivalently the morphism $F$ lead to a simplicial map $\W
\smp{H}_\bullet \to \W\smp{G}_\bullet$, the lift is only a matter
of composing $\eta$ with said map.  We prefer to present a direct
approach, which will be useful here and elsewhere in this text.

Let $V_\bullet$ be a hypercover as above, and let $\eta =
(y,h)\colon V_\bullet \to \W \smp{H}_\bullet$ be a 1-cocycle,
with $y\colon V_1\to H_0$ and $h\colon V_2\to H_1$.  Since
$\pi\colon E\to H_0$ is a sheaf epimorphism, there will be a
local lift of $y$ to $E$, namely a (generalized) cover $p_1\colon
U\to V_1$ and $e\colon U\to E$ such that
\begin{equation*}
  \xymatrix{%
    U \ar[r]^e \ar[d]_{p_1} & E \ar[d]^\pi \\
    V_1 \ar[r]_y & H_0
  }
\end{equation*}
commutes. Using \cite[V, Théorème 7.3.2]{MR0354653}, there is a
hypercover $V'_\bullet$ dominating $V_\bullet$, with a
factorization $V'_1\to U\to V_1$.  All objects will be considered
relative to $V'_\bullet$ by pull-back along the latter map. In
particular, $\eta=(y,h)$ can now be considered as a 1-cocycle
relative to $V'_\bullet$ via $V'_\bullet\to V_\bullet\to
\W\smp{H}_\bullet$.

The explicit form of the cocycle condition on $(y,h)$, the
relation $\del_H=\pi \circ \kappa$, and the injectivity of
$\imath\colon G_1\to E$ show that there must exist $g\colon
V'_2\to G_1$ such that
\begin{equation}
  \label{eq:16}
  d_1^*e = d_2^*e\, d_0^*e\, \kappa (h)\, \imath (g).
\end{equation}
Set $x = \jmath\circ e\colon V'_1\to G_0$.  We show that the pair
$(x,g)$ determines a 1-cocycle $\xi\colon V'_\bullet \to
\W\smp{G}_\bullet$.

Applying $\jmath$ to the previous relation gives the first
cocycle condition~\eqref{eq:9}.  After a pull-back to $V'_3$, and
using \eqref{eq:16} to reduce $(d_2d_3)^*e (d_0d_3)^*e\,
(d_0d_1)^*e$ in both possible ways, by a routine calculation we
obtain the equality
\begin{equation}
  \label{eq:17}
  \kappa (d_2^*h \, d_0^*h)\, \imath (d_2^*g \, d_0^*g)
  =  \kappa ( (d_3^*h)^{(d_0d_1)^*b}\, d_1^*h )\,
  \imath ((d_3^*g)^{(d_0d_1)^*x} \, d_1^*g ),
\end{equation}
so that the second cocycle condition~\eqref{eq:10} for $(x,g)$
also holds.  (This uses the fact that $\imath$ is injective and
that its image commutes with that of $\kappa$.)
\begin{remark}
  \label{rem:3}
  From~\eqref{eq:16} and~\eqref{eq:17}, it follows that
  $\tilde\eta=(e, (h,g))$ defines a 1-cocycle with values in the
  crossed module $(\kappa,\imath)\colon H_1\times G_1\to E$.
\end{remark}
\begin{remark}
  \label{rem:4}
  The technique adopted in this section can also be used to
  describe the explicit lift of a 0-cocycle with values in
  $\W\smp{H}_\bullet$ of the type discussed in
  \citepalias{ButterfliesI}.  It is an exercise to show that the
  geometric view in terms of torsors given there reduces to this
  one when trivializations are chosen.  This view is implicit in
  the proof of Theorem~\ref{thm:1} given in \citepalias[Theorem
  4.3.1]{ButterfliesI}.
\end{remark}

\subsection{Computing the map \protect {$F_*$}}
\label{sec:computing-f_colon-h1}

When $\grh \isoto [H_1\to H_0]\sptilde$, $\grg\isoto [G_1\to
G_0]\sptilde$, and $(F,\lambda)$ is expressed through the
butterfly $[H_\bullet,E,G_\bullet]$, the image of a 1-cocycle
$(y,h)$ with values in $\grh$ can be explicitly computed.  Most
of the necessary calculations follow in a straightforward way
from the explicit treatment of the equivalence between the
morphism $F$ and the butterfly provided in \citepalias[Theorem
4.3.1]{ButterfliesI} (recalled here as Theorem~\ref{thm:1}).

Recall that we have the equivalence $\grh\iso \tors (H_1,H_0)$,
and therefore, if the object $y$ corresponds to the
$(H_1,H_0)$-torsor $(Q,t)$, then $F(y)$ can be computed as
\begin{equation*}
  F(Q,t) = \shHom_{H_1}(Q,E)_t\,,
\end{equation*}
as shown in Part I.  The right-hand side is the $G_1$-torsor of
local $H_1$-equivariant lifts of $t\colon Q\to H_0$ to $E$.  In
fact it is a $(G_1,G_0)$-torsor: the section
\begin{equation*}
  s\colon \shHom_{H_1}(Q,E)_t \lto G_0
\end{equation*}
is simply the map sending a local lift $e$ of $t$ to $\jmath\circ
e$.  The morphism $h$ is the isomorphism of torsors
\begin{equation*}
  h\colon d_1^*(Q,t) \lisoto (d_2^*Q\cprod{H_1}d_0^*Q, d_2^*t d_0^*t),
\end{equation*}
so that the composite $\lambda\circ F(h)$ arises, again as
explained in Part I, from the isomorphism of $G_1$-torsors
\begin{equation*}
  \shHom_{H_1}(d_2^*Q\cprod{H_1} d_0^*Q,E)_t \lisoto
  \shHom_{H_1} (d_2^*Q,E)_t \cprod{G_1}
  \shHom_{H_1}(d_0^*Q,E)_t.
\end{equation*}

Assume the hypercover $U_\bullet$ with respect to which $(y,h)$
is defined is such that the underlying $H_1$-torsor $Q$ is
trivial, and the whole cocycle can be expressed via a 1-cocycle
with values in the crossed module $H_1\to H_0$.  Let us keep the
notation $(y,h)$ for the latter, so that now $y\in H_0(U_1)$ and
$h \in H_1 (U_2)$.

Recalling that $y\in H_0(U_1)$ corresponds to the object
$(H_1,y)$ of $\grh (U_1)$, its image under $F$ is given by:

\begin{equation}
  \label{eq:18}
  \begin{aligned}
    \shHom_{H_1}(H_1,E)_y & \lisoto E_y \\
    e & \longmapsto e(1)
  \end{aligned}
\end{equation}
where the $G_1$-torsor on the right-hand side is the ``fiber'' of
$E\to H_0$ above $y$.  It follows that the resulting cocycle with
values in $\grg$ is given by the datum of $E_y$ plus the morphism
\begin{equation}
  \label{eq:19}
  \gamma\colon E_{d_1^*y} \lisoto E_{d_2^*y} \cprod{G_1} E_{d_0^*y}.
\end{equation}
arising from the application of $(F,\lambda)$ to the first
relation in the 1-cocycle condition, i.e.
\begin{equation*}
  d_1^*y = d_2^*y\,d_0^*y\,\del h,
\end{equation*}
which really is the morphism
\begin{equation*}
  h\colon (H_1,d_1^*y) \lto (H_1, d_2^*y\,d_0^*y).
\end{equation*}
So~\eqref{eq:19} is the result of the composition
\begin{equation}
  \label{eq:20}
  E_{d_1^*y} \lto E_{d_2^*y\,d_0^*y} \lto E_{d_2^*y} \cprod{G_1} E_{d_0^*y}.
\end{equation}
A trivialization of the $G_1$-torsor $E_y$ will produce a
1-cocycle with values in the crossed module $G_1\to G_0$. More
precisely, we have:
\begin{proposition}
  \label{prop:2}
  The choice of a trivialization $e\in E_y$ amounts to a lift of
  the 1-cocycle $\eta=(y,h)\colon U_\bullet\to \W\smp{H}_\bullet$
  along the butterfly $[H_\bullet,E,G_\bullet]$, as described in
  section~\ref{sec:lift-cocycle-along}.
\end{proposition}
\begin{proof}
  One needs to show that the choice of a trivialization $e\in
  E_y$ leads to formulas~\eqref{eq:16} and~\eqref{eq:17}.
  Indeed, after pullback the choice of $e\in E_y$ yields
  $d_1^*e$, $d_2^*e$, and $d_0^*e$.

  The first morphism of~\eqref{eq:20} sends $d_1^*e$ to
  $(d_1^*e)\, \kappa (h)^{-1}$.  This is a consequence of the
  following observation: suppose we have $y = y' \, \del h$, for
  $y,y'\in H_0$ and $h\in H_1$.  Consider the diagram
  \begin{equation*}
    \xymatrix{%
      \shHom_{H_1}(H_1,E)_y \ar[d] \ar[r] &
      \shHom_{H_1}(H_1,E)_{y'} \ar[d] \\
      E_{y} \ar[r] & E_{y'}
    }
  \end{equation*}
  where the top horizontal arrow sends a local lift $e$ to
  $e\circ h^{-1}$.  Then, using~\eqref{eq:18} for the vertical
  arrows, we can calculate the bottom horizontal arrow and find
  that a section $e$ is sent by to $e\;\kappa{(h)}^{-1}$.

  Returning to the problem at hand, since the product $d_2^*e \,
  d_0^*e$ provides a trivialization of $E_{d_2^*y} \cprod{G_1}
  E_{d_0^*y}$, there must exist a $g\in G_1$ such that
  \begin{equation*}
    (d_1^*e)\; \kappa (h)^{-1} = d_2^*e \; d_0^*e \; \imath (g),
  \end{equation*}
  which clearly is the same as~\eqref{eq:16}, as wanted.

  Relation~\eqref{eq:17} follows from this last one by direct
  calculation. Alternatively, one can show that it follows from
  the cocycle condition~\eqref{eq:13} applied to the
  morphism~\eqref{eq:19}, by pulling back to $U_3$ and moving
  from $(d_1d_2)^*e$ to the product $(d_2d_3)^*e \; (d_0d_3)^*e
  \; (d_0d_1)^*e$ in the two possible ways.  The second approach
  subsumes the second.  In any event, both are straightforward
  and left to the reader.
\end{proof}

\section{Gerbes bound by a crossed module}
\label{sec:gerbes-bound-crmod}

\subsection{Recollections on gerbes}
\label{sec:recollections-gerbes}

For gerbes, our main references will be
\cite{MR49:8992,MR95m:18006}. Recall that a gerbe $\gP$ over $\s$
is by definition a stack in groupoids over $\s$ which is
``locally non-empty'' and ``locally connected.'' Following
\cite{MR1771927}, this can be expressed as follows. Let $X$ be a
``space,'' \ie a sheaf of sets, over $\s$.  A gerbe over $X$ is a
stack in groupoids $\gP$ over $\s$ equipped with a morphism
$p\colon \gP\to X$ such that both $p$ and the diagonal $\Delta
\colon \gP \to \gP\times_{X}\gP$ are (stack) epimorphisms.  The
usual definition of gerbe over $\s$ without reference to another
space is recovered by setting $X=*$. Any stack $\stack{X}$ is
equipped with a canonical morphism
\begin{equation*}
  \stack{X} \lto \pi_0(\stack{X})
\end{equation*}
which makes $\stack{X}$ into a gerbe over $\pi_0(\stack{X})$
(\cite[\S 3.19]{MR1771927} and \cite[\S 7.1]{MR95m:18006}). This
construction and its analog for 2-stacks were applied at
different points in Part I.

If $U\to *$ is a generalized cover and $G$ is a sheaf of groups
over $\s/U$, then $\gP$ is a $G$-gerbe if there exists an object
$x\in \Ob (\gP_U)$ and an isomorphism
\begin{equation*}
  G\lto \shAut_U(x)\,.
\end{equation*}
(The choice of the isomorphism is called a labeling of $\gP$ in
\cite{MR95m:18006}). It is well known from \loccit that a
$G$-gerbe gives rise to a non-abelian cohomology class with
values in the crossed module $[\iota\colon G\to \Aut
(G)]$. Essentially identical cohomology classes are shown in
\cite{MR92m:18019} to arise from $\grg$-torsors, where $\grg =
[G\to \Aut (G)]\sptilde$ is the associated gr-stack. In fact, it
is also shown in \loccit that there is an equivalence (of
2-gerbes) between $\grg$-torsors and $G$-gerbes. This section is
devoted to tie together these strands for a general crossed
module $G_1\to G_0$ of $\T$.

\subsection{Gerbes bound by a crossed module}
\label{sec:gerbes-bound-crossed}

Let $G_\bullet \colon G_1 \overset{\del}{\to} G_0$ be a crossed
module of $\T$.  The concept of gerbe bound by $G_\bullet$ is a
sort of rigidification, due to Debremaeker \cite{MR0480515}, of
the idea of $G$-gerbe recalled above.
\begin{definition}
  \label{def:3}
  A gerbe $\gP$ bound by $G_\bullet$, or equivalently, a
  $(G_1,G_0)$-gerbe, is a gerbe $\gP$ over $\s$ equipped with the
  following data:
  \begin{enumerate}
  \item a functor $\mu \colon \gP\to \tors (G_0)$;
  \item for each object $x$ of $\gP$ a functorial isomorphism $\jmath_x
    \colon \shAut (x) \isoto \mu (x) \cprod{G_0}G_1$ such that
    the diagram
    \begin{equation}
      \label{eq:21}
      \vcenter{\xymatrix@C+1pc{%
          \shAut (x) \ar[r] \ar[d]_{\jmath_x} & \shAut (\mu{(x)})
          \ar[d]^\iso \\
          \mu (x) \cprod{G_0}G_1 \ar[r]^{\id \wedge \del} &
          \mu (x) \cprod{G_0}G_0
        }}
    \end{equation}
    commutes. The right vertical morphism is
    the standard one identifying the automorphism group of a
    $G$-torsor $P$ with the twisted adjoint group
    $\operatorname{Ad} P = P\cprod{G}G$.
  \end{enumerate}
\end{definition}
Let us explicitly remark that the functoriality requirement made right
above diagram~\eqref{eq:21} means we must have, for each morphism
$f\colon x\to y$ in $\gP$, over (say) $U$, a commutative diagram
\begin{equation}
  \label{eq:22}
  \vcenter{\xymatrix@C+2pc{%
      \shAut (x) \ar[r]^{f_*} \ar[d]_{\jmath_x} & \shAut (y)
      \ar[d]^{\jmath_y} \\
      \mu (x) \cprod{G_0}G_1 \ar[r]_{\mu{(f)} \wedge \id} &
      \mu (y) \cprod{G_0}G_1
    }}
\end{equation}
where $f_*$ is defined, as usual, by sending a section $\gamma$ of
$\shAut (x)$ to $f\circ \gamma\circ f^{-1}$.  Furthermore, the obvious
cube built from~\eqref{eq:21} and~\eqref{eq:22} should commute.
\begin{example}
  \label{ex:1}
  $\tors (G_1)$ is evidently a $(G_1,G_0)$-gerbe with $\mu=\del_*$
  and $\jmath$ given by
  \begin{equation*}
    \jmath_P\colon
    P\cprod{G_1}G_1 \lisoto \del_*(P)\cprod{G_0}G_1
  \end{equation*}
  for a $G_1$-torsor $P$. $\tors (G_1)$ will be called the
  trivial $(G_1,G_0)$-gerbe when equipped with the structure just
  described.  We shall see shortly, in sect.~\ref{sec:prgr_7},
  that all $(G_1,G_0)$-gerbes are locally of this type.
\end{example}
We will denote a gerbe bound by $G_\bullet$ synthetically as
$(\gP, \mu, \jmath)$.  We have morphisms and 2-morphisms of
gerbes bound by $G_\bullet$, as follows:
\begin{definition}
  \label{def:4}
  A morphism $(F,\phi)\colon (\gP, \mu, \jmath)\to (\gP', \mu',
  \jmath')$ of gerbes bound by $G_\bullet$ is given by a morphism
  $F\colon \gP\to \gP'$ of gerbes plus a 2-morphism
  \begin{equation*}
    \xymatrix{%
      \gP \ar@/_/[dr]_\mu^{}="m" \ar[r]^F
      & \gP' \ar[d]^{\mu'}_{}="m'"
      \ar@{=>}@/^0.3pc/ "m";"m'"^\phi\\
      & \tors (G_0) }
  \end{equation*}
  such that for every object $x\in \Ob (\gP)$ the following
  diagram commutes:
  \begin{equation*}
    \xymatrix{%
      \shAut (x) \ar[r]^{F_*} \ar[d]_{\jmath_x} &
      \shAut (F(x)) \ar[d]^{\jmath'_{F(x)}} \\
      \mu (x)\cprod{G_0} G_1 \ar[r]_{\phi_x} &
      \mu'(F(x)) \cprod{G_0} G_1}
  \end{equation*}
  A 2-morphism $\theta\colon (E,\epsilon)\Rightarrow (F,\phi)$ is
  a 2-morphism of gerbes $\theta\colon E\Rightarrow F$ such that
  \begin{equation*}
    \mu'*\theta \circ \epsilon = \phi.
  \end{equation*}
\end{definition}

In ref.\ \cite{MR0480515} the definition of morphism is given in
greater generality than in Definition~\ref{def:4} above, by
allowing a \emph{strict morphism} of crossed modules. Recall that
a strict morphism $f_\bullet\colon H_\bullet\to G_\bullet$ is a
commutative diagram of group objects
\begin{equation*}
  \xymatrix{%
    H_1 \ar[r]^{f_1} \ar[d]_\del & G_1\ar[d]^\del \\
    H_0 \ar[r]_{f_0} & G_0
  }
\end{equation*}
where $f_1$ is an $f_0$-equivariant map.
\begin{definition}
  \label{def:5}
  Let $(\gP,\jmath,\mu)$ be a $(G_1,G_0)$-gerbe and $(\gQ,
  \kappa, \nu)$ an $(H_1,H_0)$-gerbe. An
  \emph{$f_\bullet$-morphism} $(F,\phi)\colon \gQ\to \gP$ is the
  datum of a morphism $F\colon \gQ\to \gP$ of gerbes plus a
  2-morphism
  \begin{equation*}
    \xymatrix{%
      \gQ \ar[d]_\nu \ar[r]^F & \gP \ar[d]^\mu_{}="b" \\
      \tors (H_0) \ar[r]_{(f_0)_*}^{}="a" & \tors (G_0)
      \ar@{=>}@/^0.3pc/ "a";"b"^\phi
    }
  \end{equation*}
  such that for each object $y$ of $\gQ$ there is a functorial
  diagram
  \begin{equation*}
    \xymatrix{%
      \shAut (y) \ar[r]^{\Aut (F)} \ar[d]_{\kappa_y} &
      \shAut (F(y)) \ar[d]^{\jmath_{F(y)}} \\
      \nu (y)\cprod{H_0} H_1 \ar[r]_\omega &
      \mu (F(y))\cprod{G_0}G_1
    }
  \end{equation*}
  where $\omega$ is the composite
  \begin{equation*}
    \nu (y) \cprod{H_0}H_1 \lto
    \nu (y) \cprod{H_0}G_1 \lisoto
    (\nu (y) \cprod{H_0}G_0 )\cprod{G_0} G_1 \lto
    \mu (F(y)) \cprod{G_0} G_1.
  \end{equation*}
\end{definition}
There is an obvious generalization of the notion of 2-morphism
too. The reader can formulate the appropriate diagram.
\begin{remark}
  An abelian crossed module is simply a homomorphism of abelian
  groups of $\s$.  Gerbes bound by crossed modules in this sense
  have appeared in refs.\ \cite{math.AG/0301304} and
  \cite{doi:10.1016/j.jpaa.2007.07.020}. As it is shown in the
  latter, the notion encompasses several well-known examples such
  that of connective structure due to Brylinski and McLaughlin
  (\cite{bry:loop}) and hermitian structure due to one of the
  authors (\cite{MR2142353}).
\end{remark}

\subsection{Local description}
\label{sec:prgr_7}

We want to explicitly show that a $(G_1,G_0)$-gerbe
$(\gP,\jmath,\mu)$ is always locally equivalent to $\tors (G_1)$
with the structure described in Example~\ref{ex:1}.

First, it will be useful to carry out a few local calculations to
translate the global structure afforded by the
$(G_1,G_0)$-structure on the gerbe $\gP$ into the operations of
the crossed module $\del\colon G_1\to G_0$.  To this end,
consider diagram~\eqref{eq:22}, and assume two trivializations
$u,v$ of the $G_0$-torsors $\mu (x)$ and $\mu (y)$ are given. It
follows that $f$ determines an element $a_f\in G_0(U)$ by
\begin{equation*}
  u\lmto \mu (f)(u) = v\,a_f.
\end{equation*}
Then $\gamma$ in $\shAut (x)$ determines, via the trivialization
$u$, an element $g\in G_1 (U)$:
\begin{equation*}
  j_x (\gamma) = u\wedge g.
\end{equation*}
From diagram~\eqref{eq:22} we have that the action of $f_*$
amounts to:
\begin{equation*}
  u \wedge g \lmto v\,a \wedge g = v\wedge g^{a^{-1}_f}.
\end{equation*}
Thus, if the trivializations are fixed, the action of $f_*$ can
be identified with the automorphism of $G_1$ given by:
\begin{equation*}
  g \lmto g^{a^{-1}_f}.
\end{equation*}
If in particular $y=x$, so that $f\in \shAut (x)$ too, then
$\jmath_x (f) = u\wedge h_f$, and by~\eqref{eq:21} we must have
$a_f = \del h_f$. Since $\jmath_x(f\circ\gamma\circ f^{-1}) =
u\wedge (h_fg{h_f}^{-1})$, it immediately follows that
\begin{equation*}
  g^{\del h_f} = {h_f}^{-1}\,g\, h_f\,.
\end{equation*}
  
Returning to the question of the local structure of $\gP$, let
$x$ be the choice of an object of $\gP_U$, for a suitable $U\to
*$.  We can assume that there exists a trivialization $s$ of the
$G_0$-torsor $\mu (x)$, refining $U$ if necessary.
\begin{lemma}
  \label{lem:3}
  The pair $(x,s)$ determines an equivalence of
  $(G_1,G_0)$-gerbes
  \begin{equation*}
    (L_{x,s},\lambda_{x,s}) \colon \gP\rvert_U \lisoto \tors (G_1).
  \end{equation*}
\end{lemma}
\begin{proof}
  The underlying functor $L_{x,s}\colon \gP\rvert_U \to \tors
  (G_1\rvert U)$ is the standard one defined by the assignment
  \begin{equation*}
    y \rightsquigarrow \shHom_\gP (x,y)
  \end{equation*}
  (see \cite{MR95m:18006,MR95b:18009}). It is the choice of $s$
  that allows us to conclude that $\shHom_\gP (x,y)$ is a
  $G_1\rvert_U$-torsor.

  Let $f\colon x\to y$ be a morphism of $\gP_U$ (over some $V\to
  U$) and let $a$ be an element of $G_0$ over $V$.  The claim is
  that the required isomorphism of $G_0\rvert_U$-torsors
  \begin{equation*}
    \lambda^{-1}_{x,s}\colon \shHom_\gP (x,y) \cprod{G_1} G_0 \lto \mu (y)
  \end{equation*}
  is defined by the assignment
  \begin{equation*}
    (f, a) \lmto \mu (f) (s) \, a.
  \end{equation*}
  Indeed, let $f$ be replaced by $f\circ\gamma$, where $\gamma$
  is an automorphism of $x$. Then there is an element $g$ of
  $G_1$ such that $\jmath_x (\gamma) = s \wedge g$, and by
  definition we have
  \begin{equation*}
    \mu (\gamma) (s) = s\, \del (g),
  \end{equation*}
  so that
  \begin{equation*}
    \mu (f\circ \gamma) (s) \, a = \mu (f) (s)\, \del (g)\,a.
  \end{equation*}
  Thus the pairs $(f\circ \gamma , a)$ and $(f, \del(g)\, a)$ map
  to the same point of $\mu (y)$, hence the claim.
\end{proof}
\begin{remark}
  For a $(G_1,G_0)$-gerbe $\gP$ choosing an object $x$ and an
  appropriate trivialization of the resulting $G_0$-torsor $\mu
  (x)$ shows that $\gP$ is in particular a $G_1$-gerbe.
\end{remark}

\subsection{The class of a gerbe bound by a crossed
  module}
\label{sec:cohom-class-gerbe}

For gerbes bound by $G_1\to G_0$ there is an analogous statement
to Proposition~\ref{prop:1}.
\begin{proposition}
  \label{prop:3}
  The elements of the pointed set $\H^1(*,G_1\to G_0)$ are in
  bijective correspondence with equivalence classes of
  $(G_1,G_0)$-gerbes over $\s$.
\end{proposition}
\begin{proof}
  Let $V_\bullet\to *$ be a hypercover such that we can choose an
  object $x\in \Ob\gP_{V_0}$ and a morphism $f \colon d_0^*x \to
  d_1^*x$ in $\Mor \gP_{V_1}$.  The choice of the pair $(x,f)$ is
  a labeling of $\gP$ relative to $V_\bullet$. Let us temporarily
  put $G=\shAut(x)$. The computations in \cite[\S
  5.2]{MR95b:18009}, show that there exists an element $\gamma$
  of $\shAut \bigl( (d_0d_1)^*x)\bigr)\iso (d_0d_1)^*G$ over
  $V_2$, defined by the diagram
  \begin{equation}
    \label{eq:23}
    \vcenter{%
      \xymatrix@+1pc{%
        (d_0d_2)^*x \ar[r]^{d_2^*f} & (d_1d_2)^*x \\
        (d_0d_1)^*x \ar[u]^{d_0^*f} &
        (d_0d_1)^*x \ar[u]_{d_1^*f} \ar[l]_\gamma
      }}
  \end{equation}
  such that the non-abelian cocycle condition holds:
  \begin{equation}
    \label{eq:24}
    \begin{aligned}
      (d_1^*f)_* & = (d_2^*f)_*\circ (d_0^*f)_* \circ
      (\iota_\gamma)   \\
      d_0^*\gamma \circ d_2^*\gamma & = \bigl( d_3^*\gamma
      \bigr)^{(d_0d_1)^*f} \circ d_1^*\gamma ,
    \end{aligned}
  \end{equation}
  where $\iota_\gamma$ denotes the image of $\gamma\in G$ in
  $\Aut (G)$ and $\gamma^f$ is a short-hand for
  $(f^{-1})_*(\gamma)$. The first equation holds over $V_2$,
  whereas the second over $V_3$.

  We can assume the $G_0\rvert_{V_0}$-torsor $P\eqdef \mu (x)$ is
  trivial over some $W\to V_0$, via some choice of $s\colon W\to
  P$. Using \cite[V, Théorème 7.3.2]{MR0354653}, we can work with
  a new hypercover $V'_\bullet$ equipped with a map $V'_\bullet
  \to V_\bullet$ such that for $n=0$ we have a factorization
  $V'_0\to W\to V_0$. Let us from now on relabel $V'_\bullet$ to
  $V_\bullet$.

  Given the foregoing assumptions, it now follows that $G=\shAut
  (x) \iso G_1\rvert_{V_0}$ and $f$ determines an element $a$ of
  $G_0$ over $V_1$, whereas $\gamma$ corresponds to an element
  $g$ of $G_1$ over $V_2$.  The local calculations of
  section~\ref{sec:prgr_7} show that~\eqref{eq:24} becomes
  \begin{equation}
    \label{eq:25}
    \begin{aligned}
      (d_1^*a) & = d_2^*a\, d_0^*a \, \del g \\
      d_0^*g \, d_2^*g & = \bigl( d_3^*g \bigr)^{(d_0d_1)^*a} \,
      d_1^*g,
    \end{aligned}
  \end{equation}
  where this time $g^a$ denotes the action of $G_0$ on $G_1$ in
  the crossed module.  This is a 1-cocycle in the same sense as
  put forward in sect.~\ref{sec:cocycles},
  equations~\eqref{eq:8}.

  The choice of a different labeling $(y,f')$, which for
  simplicity we assume to be relative to the same hypercover
  $V_\bullet$, will determine another pair $(a',g')$ satisfying
  the same non-abelian cocycle condition~\eqref{eq:25}. (To
  obtain it, we must assume as well that the $G_0$-torsor $\mu
  (y)$ is trivialized by an appropriate choice, possibly changing
  the cover again in the process.) Following \cite[\S
  5.3]{MR95b:18009} we may also assume, up to further refining
  $V_\bullet$, that we have chosen a morphism
  \begin{equation*}
    \chi \colon y \lto x
  \end{equation*}
  over $V_0$. Such choices determine an element $\eta_\chi$ of
  $\shAut (d_0^*y)$ via
  \begin{equation*}
    d_1^*\chi^{-1}\circ f\circ d_0^*\chi = f'\circ \eta_\chi.
  \end{equation*}
  Again, the calculations of section~\ref{sec:prgr_7} show that
  the pair $(\chi,\eta_\chi)$ determines, via the chosen
  trivializations, a pair $(u,h)$, with $u\in G_0(V_0)$ and $h\in
  G_1 (V_1)$. Combining the latter relation with the primed and
  unprimed versions of~\eqref{eq:24}, and using~\eqref{eq:25}, we
  arrive at the relation
  \begin{equation*}
    \begin{aligned}
      a\, d_0^*u  & = d_1^*u\, a'\, \del h \\
      g' \, (d_2^*h)^{d_0^*a'} d_0^*h & = d_1^*h\,
      g^{(d_0d_1)^*u}.
    \end{aligned}
  \end{equation*}
  By comparison with~\eqref{eq:11}, the pair $(\chi, \eta_\chi)$
  (or equivalently $(u,h)$) determines a homotopy between the two
  1-cocycles corresponding to the two different labelings of
  $\gP$.

  The quickest way to reverse the process and to reconstruct a
  $(G_1,G_0)$-gerbe starting from the datum of $(a,g)$
  satisfying~\eqref{eq:25}, relative to $V_\bullet$, is to follow
  the procedure outlined at the end of \cite[\S
  5.2]{MR95b:18009}. Briefly, from $a$ we can define a trivial
  $(G_1,G_0)$-torsor $E$ over $V_1$. Now, as observed in
  \cite{MR92m:18019} and \citepalias{ButterfliesI}, a
  $(G_1,G_0)$-torsor is in particular a $G_1$-bitorsor, hence
  refs. \cite{MR95b:18009,MR92m:18019} may be followed to descend
  $E$ (if necessary) to $V_0\times V_0$ and then to
  use~\eqref{eq:25} to conclude that $E$ defines a ``bitorsor
  cocycle'' relative to the \Cech cover $\cosk_0V_\bullet$,
  analogously to the cocycle that appeared in the proof of
  Proposition~\ref{prop:1}. From there, we can construct a
  $(G_1,G_0)$-gerbe by gluing local copies of $\tors
  (G_1\rvert_{V_0})$, considered as $(G_1,G_0)$-gerbes, according
  to Example~\ref{ex:1}. (For the gluing we must invoke the
  effectiveness of 2-descent data for $(G_1,G_0)$-gerbes.)
\end{proof}
\begin{remark}
  \label{rem:5}
  Equations~\eqref{eq:24} and~\eqref{eq:25} in the previous proof
  exhibit the same triangular and tetrahedral structure as
  equations~\eqref{eq:8} which was made explicit in
  Remark~\ref{rem:1}.  After having covered the arguments in the
  previous proof, as well as those in the one for
  Proposition~\ref{prop:1}, the tetrahedral diagrams in
  Remark~\ref{rem:1} should now appear natural.  In particular, the
  labeling for the vertices reflects the choice of a trivializing
  object and its subsequents pullbacks along the face maps of the
  hypercover.
\end{remark}
\begin{remark}
  \label{rem:6}
  Embedded in the proof of the previous proposition is the fact
  that, given two objects $x,y\in \Ob\gP_U$ above $U\in \Ob \s$,
  with chosen trivializations of the $G_0$-torsors $\mu (x)$ and
  $\mu(y)$, the $(\shAut (x), \shAut (y))$-bitorsor
  \begin{equation*}
    E_{x,y}\eqdef \shHom_{\gP}(y,x)
  \end{equation*}
  is in fact a $(G_1,G_0)$-torsor. This follows at once from the
  calculations of section~\ref{sec:prgr_7}.  From this point of
  view an arrow $f \colon y\to x$ defined over a (generalized)
  cover $V\to U$ is to be considered as a local section of such
  torsor. In particular, the assignment defined in
  section~\ref{sec:prgr_7} of $a_f\in G_0(V)$ to $f$ ought to be
  seen as the $G_1$-equivariant map
  \begin{equation*}
    s\colon E\lto G_0
  \end{equation*}
  which is part of the definition of $(G_1,G_0)$-torsor. Indeed,
  if $f$ is replaced by $f\circ \gamma$, where $\gamma\in \shAut
  (y)(V)$ and $\gamma$ is then identified with an element $g\in
  G_1(V)$, then we have
  \begin{equation*}
    a_{f\circ \gamma} = a_f\, \del g.
  \end{equation*}
\end{remark}

\subsection{Bitorsor cocycle associated to a labeling}
\label{rem:7}

According to ref.\ \cite{MR95b:18009} and Remark~\ref{rem:6}, the
proof of Proposition~\ref{prop:3} can be reformulated in terms of
the bitorsor cocycles introduced in
section~\ref{sec:bitorsor-cocycles}.  Indeed, the local
equivalence of $(G_1,G_0)$-gerbes provided by a labeling,
analyzed in section~\ref{sec:prgr_7}, in particular in
Lemma~\ref{lem:3}, determines a bitorsor cocycle as follows.  Let
\begin{equation*}
  \phi_U\colon \tors (G_1\rvert_U) \lto \gP\rvert_U
\end{equation*}
be such an equivalence, where $\gP$ is a $(G_1,G_0)$-gerbe.  Now,
let $U$ be the degree zero stage of a (generalized) cover
$U_\bullet$, and consider the two possible pull-backs $d_0^*\phi$
and $d_1^*\phi$ to $U_1$.  We obtain in this way a commutative
diagram
\begin{equation*}
  \xymatrix@C-1pc{%
    \tors (G_1\rvert_{U_1}) \ar[rr]^{\eta} \ar@/_/[dr]_{d_1^*\phi}
    && \tors (G_1\rvert_{U_1}) \ar @/^/[dl]^{d_0^*\phi}\\
    & \gP\rvert_{U_1}
  }
\end{equation*}
of $(G_1,G_0)$-gerbes which commutes up to natural isomorphism.  By
Morita theory (see \cite{MR95m:18006,MR2183393}) $\eta$ is induced by
a $G_1$-bitorsor $E$.  It is relatively easy to see that $E$ is in
fact an object of $\grg_{U_1}$, that is a $(G_1,G_0)$-torsor over
$U_1$.  The formal argument will constitute the proof of
Lemma~\ref{lem:4} below. The pull back to $U_2$ determines a
2-morphism
\begin{equation*}
  \gamma\colon
  d_1^*\eta \Rightarrow d_2^*\eta\circ d_0^*\eta \colon
  \tors(G_1\rvert_{U_2}) \lto \tors (G_1\rvert_{U_2}),
\end{equation*}
which results in the morphism of bitorsors
\begin{equation*}
  \gamma\colon d_1^*E \lto d_2^*E\cprod{G_1}d_0^*E,
\end{equation*}
with $\gamma$ to satisfy the appropriate coherence conditions
over $U_3$.  From Lemma~\ref{lem:1}, or rather its proof, we can
once again extract from $(g,\gamma)$ a cocycle with values in the
crossed module $G_1\to G_0$.

\subsection{Gerbes vs.\ torsors}
\label{sec:gerbes_torsors}

Let $\grg$ be the gr-stack $\tors (G_1,G_0)$.
Propositions~\ref{prop:1} and~\ref{prop:3} hold that
$\grg$-torsors and $(G_1,G_0)$-gerbes give rise to the same
equivalence classes of objects, in other words they are both
classified by the non-abelian cohomology set $\H^1(*,G_1\to
G_0)$.  The following is the analog of \cite[Proposition
7.3]{MR92m:18019} and the non-abelian counterpart of
\cite[Theorem 5.4.4]{doi:10.1016/j.jpaa.2007.07.020}. For the
statement, recall that $\equ$ denotes the stack of equivalences,
as defined in \cite[IV Proposition 5.2.5]{MR49:8992}.
\begin{proposition}
  \label{prop:4}
  There is a pair of quasi-inverse Cartesian 2-functors
  \begin{gather*}
    \Phi \colon \tors (\grg) \lto \gerbes (G_1,G_0), \quad
    \stX \lmto \tors (G_1) \cprod{\grg} \stX^o \\
    \intertext{and}
    \Psi \colon \gerbes (G_1,G_0)\lto \tors (\grg), \quad \gP
    \lmto \equ (\tors (G_1), \gP)
  \end{gather*}
  where for a right-$\grg$-torsor $\stX$ the symbol $\stX^o$
  denotes the opposite (left) torsor, which define a
  2-equivalence between the 2-stacks $\tors (\grg)$ and $\gerbes
  (G_1,G_0)$ over $\s$.
\end{proposition}
In fact the pair defines a 2-equivalence between neutral 2-gerbes
over $\s$.  For the proof, the following lemma, which is also of
independent interest, is needed:
\begin{lemma}
  \label{lem:4}
  There is an equivalence of gr-stacks
  \begin{equation*}
    \grg \lisoto \equ\, (\tors (G_1), \tors(G_1))
  \end{equation*}
  where $\tors (G_1)$ is considered as a $(G_1,G_0)$-gerbe in the
  manner described by Example~\ref{ex:1}.
\end{lemma}
\begin{proof}
  The functor in the statement is the one sending the
  $(G_1,G_0)$-torsor $(E,s)$ to the equivalence
  \begin{equation*}
    P \lmto P\cprod{G_1} E
  \end{equation*}
  where, according to \cite{MR92m:18019}, recalled in
  \citepalias[\S 3.4.8]{ButterfliesI}, $E$ is a $G_1$-bitorsor
  using the left $G_1$-action defined as $g\cdot e = e
  g^{s(e)}$. The functor is clearly fully faithful.
  
  Let $(F,\phi)\colon \tors (G_1)\to \tors (G_1)$ be an
  equivalence of $(G_1,G_0)$-gerbes (see Definition~\ref{def:4}).
  Recall that by standard arguments of Morita theory, the
  underlying functor $F$ determines and is determined, up to
  equivalence, by a $G_1$-bitorsor $E$ so that for any right
  $G_1$-torsor $P$ there is an isomorphism
  \begin{equation*}
    F(P) \iso P\cprod{G_1} E.
  \end{equation*}
  $E$ is simply the image under $F$ of the trivial torsor $G_1$.
  By Definition~\ref{def:4}, this must be compatible with
  $\del_*\colon \tors (G_1)\to \tors(G_0)$, so there must exist
  an isomorphism
  \begin{equation*}
    \phi_P \colon P\cprod{G_1}G_0 \lisoto
    \bigl(P\cprod{G_1} E\bigr) \cprod{G_1} G_0
  \end{equation*}
  for all torsors $P$. If in particular $P=G_1$, it reduces to
  \begin{equation*}
    \phi_{G_1} \colon G_0 \lisoto E\cprod{G_1} G_0,
  \end{equation*}
  that is $E$ must be equipped, as a right $G_1$-torsor, with a
  trivialization of its extension to a $G_0$-torsor. Thus $E$ is
  a $(G_1,G_0)$-torsor, and it is relatively easy to verify that
  the resulting left $G_1$-torsor structure recalled above is the
  same as the original one.
\end{proof}
\begin{proof}[Main lines of the proof of
  Proposition~\ref{prop:4}]
  The proof closely mirrors the one in \cite[Proposition
  7.3]{MR92m:18019}, except for the details pertaining to the
  $(G_1,G_0)$-gerbe structure.
  
  By Lemma~\ref{lem:4}, $\grg$ acts on the right on $\Psi
  (\gP)$. As observed in \loccit, for any two equivalences $F,F'$
  we have $F' \iso F \circ (F^{-1}\circ F')$, for a choice
  $F^{-1}$ of the quasi-inverse to $F$, and $F^{-1}\circ F$ is an
  auto-equivalence of $\tors (G_1)$.  Furthermore, $\Psi (\gP)$
  is locally non void, since from~\ref{sec:prgr_7} the choice of
  an object $x$ of $\gP$ and of a trivialization $s$ of $\mu (x)$
  over some $U\in \Ob (\s)$ determines an equivalence $\tors
  (G_1) \isoto \gP$ of $(G_1,G_0)$-gerbes over $U$.

  As for $\Phi (\stX)$, it is a gerbe since, as already noted in
  \loccit, the very fact that $\stX$ is itself locally equivalent
  to $\grg$ shows that $\Phi (\stX)$ is locally equivalent to
  $\tors (G_1)$.

  It is to be shown that $\Phi (\stX)$ actually is a
  $(G_1,G_0)$-gerbe.  To this end, let $\mu \colon \Phi (\stX)
  \to \tors (G_0)$ be defined by
  \begin{equation}
    \label{eq:26}
    \mu (P,X) \eqdef \del_* (P) = P\cprod{G_1} G_0.
  \end{equation}
  If the triple $(\alpha, g, \beta)$, where $g=(E,s)$ denotes a
  $(G_1,G_0)$-torsor, represents a morphism
  \begin{equation*}
    (P_1,X_1) \lto (P_2,X_2)
  \end{equation*}
  in $\Phi (\stX)$ as described in~\ref{sec:prgr_10}, then
  $\del_*([\alpha,g,\beta])$ is defined to be the composition
  \begin{equation}
    \label{eq:27}
    P_1\cprod{G_1} G_0 \xrightarrow{\alpha\wedge \id_{G_0}}
    (P_2\cprod{G_1} E) \cprod{G_1} G_0 \lisoto
    P_2\cprod{G_1} (E \cprod{G_1} G_0)
    \xrightarrow{\id_{P_2}\wedge s} P_2\cprod{G_1} G_0.
  \end{equation}
  It is immediately checked that it does not depend on the
  specific choice of the triple representing the morphism.

  For two morphisms $(P_1,X_1) \to (P_2,X_2)$ and $(P_2,X_2) \lto
  (P_3,X_3)$ composed as in~\ref{sec:prgr_10}, a diagram chase,
  using Mac Lane's pentagon, reveals that the composition of the
  corresponding images~\eqref{eq:27} equals (as expected) the
  image of the composition under $\mu$.

  Having defined $\mu$, it must be proved that there is a
  functorial isomorphism
  \begin{equation}
    \label{eq:28}
    \jmath_{P,X}\colon \shAut (P,X) \lisoto \mu (P,X)\cprod{G_0} G_1,
  \end{equation}
  as per Definition~\ref{def:3}.  Note that from~\eqref{eq:26} it
  follows that:
  \begin{equation*}
    \mu (P,X)\cprod{G_0} G_1 \iso P\cprod{G_1} G_1 \iso
    \shAut (P), 
  \end{equation*}
  so that~\eqref{eq:28} amounts to showing that:
  \begin{equation*}
    \shAut (P,X) \iso \shAut(P).
  \end{equation*}
  This actually follows from the fact that the choice of the
  object $X$ of $\stX^o$ establishes a local equivalence with
  $\grg$, and hence one of $\Phi (\stX)$ with $\tors
  (G_1)$. Explicitly, and somewhat more precisely, an
  automorphism of $(P,X)$ is given by a triple $(\alpha, g,
  \beta)$ such that
  \begin{equation*}
    \alpha\colon P \lto P\cprod{G_1}E \qquad
    \beta\colon g\cdot X \lto X\,, \quad g = (E,s).
  \end{equation*}
  Since $\stX$ is a torsor, it follows there must be an arrow
  \begin{equation*}
    \gamma \colon g \lto I_\grg,
  \end{equation*}
  in $\grg$, that is the $(G_1,G_0)$-torsor $(E,s)$ is isomorphic
  to the trivial $(G_1,G_0)$-torsor $(G_1,1)$. It follows that
  the triple $(\alpha, g, \beta)$ is equivalent in the sense
  of~\ref{sec:prgr_10} to $(\alpha', I_\grg, l_X)$, where $l_X$
  is the structural functorial isomorphism
  \begin{equation*}
    l_X \colon I_\grg\cdot X \lisoto X
  \end{equation*}
  which is part of the definition of $\grg$-torsor. On the other
  hand, $\alpha'$ is the composition $(\id_P\cdot \gamma)\circ
  \alpha \colon P \to P\cprod{G_1} G_1\iso P$, which is the
  sought-after element of $\shAut (P)$. It is clear the
  requirements of Definition~\ref{def:3} and in~\ref{sec:prgr_7}
  are met.

  As a last point, since $\tors (G_1)\cprod{\grg}\stX^o$ is
  actually defined by a process of stackification, it should also
  be checked that $\mu$ as defined glues along descent data. If
  $(P,X)$ is an object defined over $V$ with a morphism
  \begin{equation*}
    \phi\colon d_0^*(P,X) \lto d_1^*(P,X)
  \end{equation*}
  over, say, $V\times_U V$ such that the cocycle condition
  \begin{equation*}
    d_1^* \phi= d_2^*\phi\circ d_0^*\phi
  \end{equation*}
  holds, the definition~\eqref{eq:26} should give rise to a
  well-defined $G_0$-torsor over $U$ (via descent in $\tors
  (G_0)$).  Writing $\phi$ as being represented by a triple
  $(\alpha, g, \beta)$, the descent datum above gives rise to two
  diagrams
  \begin{gather*}
    \xymatrix{%
      d_2^*g\cdot ( d_0^*g \cdot (d_0d_1)^*X) \ar[r] \ar[d]^\wr &
      d_2^*g \cdot (d_0d_2)^*X \ar[r] &
      (d_1d_2)^*X \\
      (d_2^*g\cdot d_0^*g )\cdot (d_0d_1)^*X \ar[rr] & &
      d_1^*g\cdot (d_0d_1)^*X \ar[u]
    } \\
    \intertext{and} \xymatrix{%
      (d_0d_1)^*P \ar[r] \ar[d] & (d_0d_2)^*P \cdot d_0^*g \ar[r]
      &
      ((d_1d_2)^*P \cdot d_2^*g)\cdot d_0^*g  \ar[d]^\wr \\
      (d_1d_2)^*P & & (d_1d_2)^*P \cdot (d_2^*g\cdot d_0^*g )
      \ar[ll] }
  \end{gather*}
  Applying $\mu$ produces an object $P\cprod{G_1}G_0$ over $V$ and a
  morphism $d_0^*P\cprod{G_1}G_0 \to d_1^*P\cprod{G_1}G_0$ of
  type~\eqref{eq:27} over $V\times_U V$. After having applied $\mu$ to
  the second diagram above, another long but totally straightforward
  diagram chase leads to a corresponding cocycle condition. Hence
  $P\cprod{G_1}G_0$ can be descended to a $G_0$-torsor over $U$, as
  wanted.
\end{proof}

Passing to classes of equivalences, we have the identifications
\begin{equation*}
  \bigl[ \tors (\grg) \bigr] \iso
  \bigl[ \gerbes (G_1,G_0) \bigr] \iso
  \H^1 (\ast, G_1\to G_0),
\end{equation*}
where $[\cdot]$ denotes taking classes of equivalences of objects
over $\ast$.  The first identification is of course induced by
$\Phi$ (and its inverse by $\Psi$).  It follows at once from
Proposition~\ref{prop:4} and from Propositions~\ref{prop:1}
and~\ref{prop:3} that the above identifications constitute a
commutative diagram, namely the isomorphism induced by $\Phi$ is
compatible with taking cohomology classes, so that the induced
map on $\H^1$ is the identity. We record this as a lemma.
\begin{lemma}
  \label{lem:5}
  The maps induced by $\Phi$ and $\Psi$ preserve equivalence
  classes.
\end{lemma}
For future use, it is nevertheless convenient to have a
computational verification.
\begin{proof}[Proof of the lemma]
  If $\stX$ is a $\grg$-torsor, then the choice of an object $x$
  in the fiber $\stX_U$ over $U$ establishes an equivalence
  \begin{equation*}
    \stX\rvert_U \lisoto \grg\rvert_U
  \end{equation*}
  which gives (see \cite{MR92m:18019} and the proof of
  Proposition~\ref{prop:4})
  \begin{equation*}
    \begin{split}
      \Phi (\stX\rvert_U) = \tors (G_1\rvert_U)
      \cprod{\grg\rvert_U} \stX^{o}\rvert_U & \lisoto \tors
      (G_1\rvert_U) \cprod{\grg\rvert_U} \grg^{o}\rvert_U \\
      & \lisoto \tors (G_1\rvert_U).
    \end{split}
\end{equation*}
  Explicitly, an inverse equivalence is given by:
  \begin{equation*}
    \begin{aligned}
      \phi_U\colon \tors (G_1\rvert_U) & \lisoto \tors
      (G_1\rvert_U)
      \cprod{\grg\rvert_U} \stX^{o}\rvert_U \\
      P & \lmto (P,x).
    \end{aligned}
  \end{equation*}
  According to section~\ref{rem:7}, this equivalence will
  determine a bitorsor cocycle for the gerbe $\Phi (\stX)$, which
  we want to identify with the one determined by the choice of
  the object $x$ of $\stX$. Indeed, let the latter be given by
  the pair $(g,\gamma)$, with $g=(E,s)$ is a $(G_1,G_0)$-torsor
  over $U=U_0$, as in the proof of Proposition~\ref{prop:1}. From
  the morphism
  \begin{equation*}
    \xi \colon d_0^*x \lisoto d_1^*x\cdot g
  \end{equation*}
  in $\stX_{U_1}$ consider the morphism ($g^*$ is a choice of the
  inverse for $g$):
  \begin{equation*}
    d_0^*x \cdot g^* \lisoto
    \bigl( d_1^*x\cdot g \bigr)\cdot g^* \lisoto
    d_1^*x\cdot \bigl(g\cdot g^*\bigr) \lisoto d_1^*x,
  \end{equation*}
  which by definition corresponds to a morphism $\xi^o$ in
  $\stX^o$:
  \begin{equation*}
    \xi^0 \colon g\cdot d_0^*x \lto d_1^*x.
  \end{equation*}
  By the definition of contracted product given in
  sect.~\ref{sec:prgr_10}, we have that the triple $(id,g,\xi^o)$
  determines a morphism
  \begin{equation*}
    d_0^*\phi \bigl(P\cprod{G_1} E\bigr)
    = (P\cprod{G_1} E, d_0^*x) \coin
    (P\cdot g, d_0^*x) \lisoto (P, d_1^*x)
    =d_1^*\phi \bigl( P\bigr).
  \end{equation*}
  By comparison with the results of section~\ref{rem:7}, we see
  that resulting self-equivalence of $\tors (G_1\rvert_{U_1})$ is
  indeed given by $g=(E,s)$, as wanted.

  In the opposite direction, let $\gP$ be a $(G_1,G_0)$-gerbe.
  If $x$ is an object of $\gP_U$, this choice will determine as
  in section~\ref{rem:7} a bitorsor cocycle $(g,\gamma)$,
  relative to some cover of $U$, where we write again $g=(E,s)$.
  In view of Lemma~\ref{lem:4}, and the definition of $\Psi$, it
  is immediate that the bitorsor cocycle for the $\grg$-torsor
  $\equ (\tors (G_1),\gP)$ (relative to the trivialization
  induced by $x$) is still $(g,\gamma)$.
\end{proof}
\begin{remark}
  The preceding proof in fact shows that both $\Phi$ and $\Psi$
  act as identities on bitorsor cocycles, thereby implying the
  statement of the lemma.
\end{remark}

\section{Extension of gerbes along a butterfly}
\label{sec:Extension-butterfly}

Functoriality of cohomology under a change of coefficients is one
of the most important properties which are required to hold in
the realm of non-abelian cohomology.  In the case of groups it is
well known that the map $\H^1(*, H) \to \H^1(*, G)$ induced by a
homomorphism $\delta\colon H\to G$ is realized by the standard
extension of torsors $\delta_*\colon \tors (H) \to \tors(G)$,
which sends an $H$-torsor $P$ to its extension $\delta_* P =
P\cprod{H}G$.  (In fact there is a $\delta$-morphism $P\to
\delta_*P$, see \cite{MR49:8992}.)

In the case of a morphism $F\colon \grh\to \grg$ of gr-stacks,
the categorification of the above extension of torsors yields the
required map $\H^1(*, \grh) \to \H^1(*, \grg)$, see ref.\
\cite{MR92m:18019}.  These matters are briefly recalled, mostly
for convenience, in sect.~\ref{sec:extension-torsors} below.
Just note that the categorification entails considering the
morphism of 2-gerbes $F_*\colon \tors(\grh) \to \tors(\grg)$
given by sending the $\grh$-torsor $\stY$ to $F_*\stY = \stY
\cprod{\grh}\grg$. In view of the equivalence between torsors and
gerbes stated in Proposition~\ref{prop:4}, this picture could be
reinterpreted in terms of gerbes bound by crossed modules, albeit
not in an immediately explicit form.

Our purpose is to remedy this by putting forward a better and
more explicit picture which rests on the equivalence (cf.\
Theorem \ref{thm:1}) between morphisms of gr-stacks and
butterflies between crossed modules, and on the interpretation of
the first non-abelian cohomology group with values in a gr-stack
as equivalence classes of gerbes.  The procedure to be expounded
below starts with a gerbe bound by the crossed module $H_\bullet$
and uses the butterfly representing $F\colon \grh\to \grg$ to
construct in a fairly explicit way a gerbe bound by $G_\bullet$,
compatibly with the categorification above.  It builds upon and
improves an earlier notion of Debremaeker~\cite{MR0480515}.

\subsection{Extension of torsors}
\label{sec:extension-torsors}

A morphism $F\colon \grh\to \grg$ of gr-stacks induces a morphism
\begin{equation*}
  F_* \colon \tors (\grh) \lto \tors (\grg)
\end{equation*}
between the corresponding 2-gerbes of torsors. The definition of
$F_*$ is the categorification of the standard ``extension of the
structural group'' for torsors, namely if $\stY$ is an
$\grh$-torsor, then we define
\begin{equation*}
  F_*(\stY) = \stY\cprod{\grh}\!\grg.
\end{equation*}
This was extensively used---without definition, but referring
instead to~\cite{MR92m:18019}---in Part I. Passing to cohomology,
that is, to isomorphism classes of objects, it is clear that
there results a corresponding maps of pointed sets:
\begin{equation*}
  \H^1(*,\grh) \lto \H^1(*,\grg).
\end{equation*}
Indeed, still according to \cite{MR92m:18019}, this is the
enabling framework to interpret the functoriality of non-abelian
cohomology with values in a crossed-module. Insofar as cohomology
only depends on the quasi-isomorphism class of the coefficient,
and every gr-stack is equivalent to one associated to a crossed
module, this covers the general case.

Let $\grh$ and $\grg$ be associated to crossed modules $H_1\to
H_0$ and $G_1\to G_0$, respectively.  In view of the equivalence
stated in Proposition~\ref{prop:4}, there is an abstract
description of $F_*$ in terms of gerbes.  Following ref.\
\cite{MR92m:18019}, let us use the notation $F_{**}$ for the
morphism $\gerbes (H_1,H_0) \to \gerbes (G_1, G_0)$ resulting
from $F_*$ via the following diagram:
\begin{equation*}
  \xymatrix{%
    \tors (\grh) \ar[d]_\Phi \ar[r]^{F_*} & \tors (\grg)
    \ar[d]^{\Phi} \\
    \gerbes (H_1,H_0) \ar[r]_{F_{**}} & \gerbes (G_1,G_0)
  }
\end{equation*}
The definition is $F_{**} = \Phi\circ F_*\circ \Psi$, and the above
diagram commutes up to natural isomorphism.

It is clear that modulo the obvious isomorphism above the
statement of Lemma~\ref{lem:5}, $F_*$ and $F_{**}$ induce the
same map $\H^1(*,\grh)\to \H^1(*,\grg)$.

Unfortunately, without additional input, $F_{**}$ cannot be
easily characterized.  If $\stY$ is again an $\grh$-torsor, a
simple manipulation gives that the gerbe $\Phi (F_*(\stY))$ is
equivalent to $\tors(G_1)\cprod{\grh}\!\stY^0$, where $\tors
(G_1)$ carries an $\grh$-action via
\begin{equation*}
  \grh \overset{F}{\lto} \grg \lisoto
  \equ (\tors (G_1),\tors (G_1)).
\end{equation*}
Thus, if $\gQ$ is an $(H_1,H_0)$-gerbe, the previous observation
suggests that its image under $F_{**}$ is
\begin{equation*}
  F_{**} (\gQ) = \tors (G_1) \cprod{\grg}\!\Psi (\gQ)^o =
  \tors (G_1) \cprod{\grh}\! \equ\, (\tors(H_1),\gQ)^o.
\end{equation*}
To improve on this picture, we propose to provide an explicit
characterization of $F_{**}$ by employing the butterfly
construction of the morphism $F\colon \grh \to \grg$.

\subsection{Debremaeker's extension along strict morphisms}
\label{sec:Extension-stric}

Let $f_\bullet\colon H_\bullet\to G_\bullet$ be a \emph{strict}
morphism of crossed modules, as in Definition~\ref{def:5}. Let
$(\gP, \jmath, \mu)$ be an $(H_1,H_0)$-gerbe. In
\cite{MR0480515}, Debremaeker proved that there exists a
$(G_1,G_0)$-gerbe $(\gP', \jmath', \mu')$ and an
$f_\bullet$-morphism $\gP \to \gP'$.

The gerbe $(\gP',\jmath',\mu')$ is constructed in two steps.
First, a fibered category $\gP^*$ is defined with the same
objects as $\gP$ and morphisms given by the extension of torsors
\begin{equation}
  \label{eq:29}
  \shHom^* (y,x) \eqdef
  \shHom_\gP (y,x) \bigwedge^{\mu (y)\cprod{H_0}H_1}
  \bigl(\mu (y)\cprod{H_0} G_1\bigr),
\end{equation}
for any two objects $x,y$ of $\gP$.  Note that in the above
formula, to define $\mu (y)\cprod{H_0}G_1$, $G_1$ is considered
as an $H_0$ object via the homomorphism $f_0\colon H_0\to G_0$,
and that the homomorphism $\id_{\mu (y)}\wedge f_1 \colon \mu
(y)\cprod{H_0}H_1 \to \mu (y)\cprod{H_0}G_1$ is used for the
extension.  Then, the second step is to define $\gP'$ as the
stack associated to $\gP^*$.  The $f_\bullet$-morphism from
$\gP\to \gP'$ is induced by the corresponding $\gP\to \gP^*$
simply given by the identity on objects and the map $f \mapsto
(f,1)$ on morphisms.

To see that $\gP'$ is a $(G_1,G_0)$-gerbe, one can argue that a
choice of trivializations of $\mu (y)$ and $\mu (x)$ above makes
$\shHom_\gP (y,x)$ into an $(H_1,H_0)$-torsor.  Consequently,
$\shHom^*(y,x)\iso\shHom_\gP (y,x)\cprod{H_1}G_1$ is a
$(G_1,G_0)$-torsor. The conclusion follows from the application
of this argument to the class of $\gP$ constructed in
Proposition~\ref{prop:3}. Still according to the proposition, the
modified cohomology class according to~\eqref{eq:29} is therefore
the class of a $(G_1,G_0)$-gerbe.

To elaborate further, according to \cite{MR0480515}, there is a
composition
\begin{equation*}
  \shHom^* (y,x) \times \shHom^* (z,y) \lto \shHom^*(z,x)
\end{equation*}
defined as follows. If $\gamma_y$ is an element of $\shAut(y)\iso
\mu(y)\cprod{H_0}G_1$, and similarly for $\gamma_z$, then the
composition law is defined as:
\begin{equation*}
  ((f,\gamma_y) , (g, \gamma_z)) \lmto
  (f\circ g, \mu(g)^{-1} (\gamma_y)\gamma_z),
\end{equation*}
where $\mu(g)^{-1}$ is a short-hand for the homomorphism of group
objects
\begin{equation*}
  \mu(y)\cprod{H_0} G_1 \lto \mu(z)\cprod{H_0}G_1
\end{equation*}
induced by $\mu(g)^{-1}\colon \mu(y)\to \mu(z)$. Note that the
functor $\mu'\colon \gP'\to \tors (G_0)$ is simply induced by the
composition of $\mu$ with
\begin{equation*}
  (f_0)_*\colon \tors (H_0) \lto \tors (G_0),
\end{equation*}
in other words to any object $x$ we assign
$\mu(x)\cprod{H_0}G_0$.  Moreover, from~\eqref{eq:29} it
immediately follows that if $y=x$ then
\begin{equation*}
  \shAut^*(x) \iso \mu(x)\cprod{H_0}G_1\iso
  \bigl(\mu (x) \cprod{H_0}G_0\bigr)\cprod{G_0} G_1,
\end{equation*}
which gives the required isomorphism $\jmath'_x$. All the
necessary requirements can be easily checked by the reader as an
exercise.

It is also not hard to realize that Debremaeker's construction is
actually functorial with respect to morphisms (and 2-morphisms)
of $(H_1,H_0)$-gerbes (see \cite{MR0480515} for details).  This
provides us with a 2-functor
\begin{equation}
  \label{eq:30}
  F^0_{+} \colon \gerbes(H_1,H_0) \lto \gerbes (G_1,G_0)
\end{equation}
which we seek to generalize in
section~\ref{sec:extension-butterfly}, to a morphism which is not
necessarily assumed to be strict.

\begin{remark}
  \label{rem:8}
  The object $E_{x,y}=\shHom_\gP(y,x)$ is a
  $(\mu(x)\cprod{H_0}H_1,\mu(y)\cprod{H_0}H_1)$-bitorsor.  It
  must be characterized (see again \cite{MR92m:18019}) by a
  $\mu(y)\cprod{H_0}H_1$-equivariant morphism
  \begin{equation*}
    E_{x,y} \lto
    \shIsom (\mu(x)\cprod{H_0}H_1,\mu(y)\cprod{H_0}H_1)
    \iso \shHom_{H_0} (\mu(y),\mu(x))
  \end{equation*}
  from $E_{x,y}$ considered as a right torsor.  This map is
  simply given by
  \begin{equation}
    \label{eq:31}
    f \lmto \mu(f)^{-1}
  \end{equation}
  where we use the same short-hand notation as above.
  Consequently, $E_{x,y}^*=\shHom^*(y,x)$ given by~\eqref{eq:29}
  has the structure of
  $(\mu(x)\cprod{H_0}G_1,\mu(y)\cprod{H_0}G_1)$-bitorsor, since
  by~\eqref{eq:31} above we get an obvious map
  \begin{equation*}
    \shIsom (\mu(x)\cprod{H_0}H_1,\mu(y)\cprod{H_0}H_1)
    \lto
    \shIsom (\mu(x)\cprod{H_0}G_1,\mu(y)\cprod{H_0}G_1),
  \end{equation*}
  which is equivariant with respect to
  \begin{equation*}
    \id \wedge f_1\colon \mu(x)\cprod{H_0}H_1\lto \mu(x)\cprod{H_0}G_1.
  \end{equation*}
  According to \cite[Proposition 2.11]{MR92m:18019}, this is what
  is required to obtain an extension of bitorsors. Thus an
  alternative way to construct the gerbe $\gP'$ is to start from
  the bitorsor cocycle $E^*$ as described in \cite{MR95b:18009}.
\end{remark}

\subsection{Extension along a butterfly}
\label{sec:extension-butterfly}

Let now $F\colon \grh\to \grg$ be a general morphism of
gr-stacks, and let $[H_\bullet, E, G_\bullet]$ be the
corresponding butterfly~\eqref{eq:6}, under the equivalence
theorem~\ref{thm:1} (we assume equivalences $\grh \iso [H_1\to
H_0]$ and $\grg\iso [G_1\to G_0]$ have been chosen).  Let also
$E_\bullet \colon H_1\times G_1\to E$ be the intermediate crossed
module, quasi-isomorphic to $H_\bullet$.  Recall
that there is a fraction
\begin{equation*}
    \xymatrix@1{%
    H_\bullet & E_\bullet \ar[l]_\sim \ar[r] & G_\bullet}
\end{equation*}
which, denoting by
$\gre$ the gr-stack associated to $E_\bullet$, factors the
morphism $F$ into
\begin{equation*}
  \grh \longleftarrow \gre \lto \grg,
\end{equation*}
where the left-pointing arrow is an equivalence.  Also, let
$(\gQ,k, \nu)$ be a gerbe bound by $H_\bullet$. The following
theorem generalizes the analogous statement of \cite[Theorem, \S
2, p. 66]{MR0480515}.
\begin{theorem}
  \label{thm:2}
  For a butterfly $[H_\bullet, E, G_\bullet]$ as above, and a
  gerbe $\gQ$ bound by $H_\bullet$, there exists a gerbe $\gP$
  bound by $G_\bullet$. The construction of $\gP$ is purely in
  terms of the butterfly $[H_\bullet,E,G_\bullet]$.
\end{theorem}
\begin{proof}
  The construction of the gerbe $\gP$ is carried out in two
  steps:
  \begin{itemize}
  \item first, construct an intermediate gerbe bound by
    $E_\bullet$;
  \item second, apply the construction of
    sect.~\ref{sec:Extension-stric} to the strict morphism
    \begin{equation}
      \label{eq:32}
      \vcenter{%
        \xymatrix{%
          H_1\times G_1 \ar[d] \ar[r]^{\mathrm{pr}_2} & G_1
          \ar[d]^\del \\
          E  \ar[r]_\jmath & G_0
        }}
    \end{equation}
    to obtain the required $(G_1,G_0)$-gerbe $\gP$.
  \end{itemize}
  To realize the first step, let us consider the gerbe:
  \begin{equation*}
    \gQ'\eqdef
    \gQ \times_{\tors (H_0)} \tors (E),
  \end{equation*}
  where the fiber product is of course taken in the sense of
  stacks: an object of $\gQ'$ is a triple $(x,f,P)$, where $x$ is
  an object of $\gQ$, $P$ is an $E$-torsor, and $f$ is an
  isomorphism
  \begin{equation*}
    f\colon \nu (x) \lisoto \pi_*(P)=P\cprod{E}H_0.
  \end{equation*}
  There is an obvious morphism
  \begin{math}
    \gQ' \lto \gQ
  \end{math}
  given by the projection to the first factor.  The proof is
  completed by showing that $\gQ'$ is bound by $E_\bullet$, which
  we state in the following lemma.
\end{proof}
\begin{lemma}
  \label{lem:6}
  The gerbe $\gQ'$ is bound by $E_\bullet\colon H_1\times G_1\to
  E$.
\end{lemma}
\begin{proof}
  Indeed, first of all there is a morphism
  \begin{equation*}
    \nu' \colon \gQ' \lto \tors (E)
  \end{equation*}
  given by the projection to the second factor, and, second,
  there is a functorial isomorphism
  \begin{equation}
    \label{eq:33}
    k' \colon \shAut (x,f, P) \lisoto P\cprod{E} (H_1\times G_1)
    \iso (P\cprod{E}H_1) \times (P\cprod{E} G_1)
  \end{equation}
  satisfying the requirements in Definition~\ref{def:3}.  To see
  this, observe that by the very definition of stack fiber
  product an automorphism of $(x,f, P)$ is given by a pair
  \begin{equation*}
    \phi\colon x \lto x\qquad \alpha\colon P\lto P
  \end{equation*}
  such that
  \begin{equation*}
    \xymatrix{%
      \nu (x) \ar[r]^f \ar[d]_{\nu (\phi)} & { P\cprod{E}H_0}
      \ar[d]^{\alpha\wedge \id} \\
      \nu (x) \ar[r]_f& { P\cprod{E}H_0}
    }
  \end{equation*}
  commutes.  In other words, $f$ determines an isomorphism
  \begin{equation*}
    f_* \colon \shAut (\nu (x)) \lisoto \shAut (P\cprod{E}H_0)
  \end{equation*}
  so that $f_*(\nu (\phi)) = \alpha\wedge \id_{H_0}$. Note that
  it coincides with
  \begin{equation*}
    f\wedge \id_{H_0}\colon \nu (x) \cprod{H_0} H_0 \lto
    \bigl( P\cprod{E}H_0 \bigr) \cprod{H_0} H_0 \iso
    P\cprod{E}H_0
  \end{equation*}
  modulo the canonical isomorphism which identifies, for any
  $G$-torsor $R$, $\shAut (R)$ with $R\cprod{G}G$.  Thus, the
  following diagram
  \begin{equation*}
    \xymatrix{%
      \shAut (x) \ar[r]^{k_x} \ar[d]_\nu &
      { \nu (x)\cprod{H_0}H_1}
      \ar[r]^(0.4){f\wedge \id}
      \ar[d]_{\id\wedge \del} &
      { \bigl( P\cprod{E}H_0 \bigr)
        \cprod{H_0}H_1} 
      \ar[r]^(0.6)\iso
      \ar[d]^{\id\wedge\del} &
      { P\cprod{E}H_1} \ar[d]^{\id\wedge\del}\\
      \shAut (\nu (x)) \ar[r]^\iso &
      { \nu (x)\cprod{H_0}H_0}
      \ar[r]^(0.4){f\wedge \id} &
      { \bigl( P\cprod{E}H_0 \bigr)
        \cprod{H_0}H_0} 
      \ar[r]^(0.6)\iso &
      P\cprod{E}H_0
    }
  \end{equation*}
  commutes. It shows that there is an isomorphism
  \begin{equation}
    \label{eq:34}
    \shAut (x,f,P) \lisoto
    P\cprod{E}H_1 \times_{(P\cprod{E}H_0)} P\cprod{E}E
    \iso
    P\cprod{E} \bigl( H_1\times_{H_0}E \bigr),
  \end{equation}
  and moreover, everything is clearly functorial.  From the
  butterfly~\eqref{eq:6} it readily follows that
  \begin{equation*}
    H_1\times_{H_0}E \iso H_1\times G_1,
  \end{equation*}
  so that \eqref{eq:34} is the promised
  isomorphism~\eqref{eq:33}, and this concludes the proof of the
  lemma.
\end{proof}
\begin{remark}
  \label{rem:9}
  Since the strict morphism~\eqref{eq:32} involves just the
  projection from $H_1\times G_1$ to $G_1$, the effect
  of~\eqref{eq:29} is to just kill off the $H_1$-part of the
  automorphisms.  More precisely, given two objects $(x,f,P)$ and
  $(y,g,Q)$ of $\gQ'$, the torsor
  \begin{equation*}
    \shHom_{\gQ'} \bigl((y,g,Q),(x,f,P)\bigr)
  \end{equation*}
  is isomorphic, via \eqref{eq:33}, to a product. In this simpler
  situation, the net effect of~\eqref{eq:29} is that of killing
  the factor relative to $P\cprod{E}H_1$.
\end{remark}
\begin{remark}
  \label{rem:10}
  The construction of the gerbe $\gP$ provided by
  Theorem~\ref{thm:2} can be described by the diagram
  \begin{equation*}
    \xymatrix@-1pc{%
      & \gQ' \ar@/_0.2pc/[dl] \ar@/^0.2pc/[dr] \\
      \gQ && \gP
    }
  \end{equation*}
\end{remark}
Both steps in the construction of the gerbe $\gP$ in the proof of
Theorem~\ref{thm:2} are (2-)functorial: this is clear for the
first step involving the fiber product construction of the gerbe
\begin{equation*}
  \gQ' = \gQ\times_{\tors (H_0)} \tors (E)
\end{equation*}
bound by $E_\bullet$, and for the second step it follows from the
functoriality of Debremaeker's construction itself, recalled in
sect.~\ref{sec:Extension-stric}.

Let $F\colon \grh\to \grg$ be the morphism of gr-stacks
corresponding to the butterfly $[H_\bullet, E, G_\bullet]$. By
the above, we have another 2-functor. We state it as follows:
\begin{definition}
  \label{def:6}
  Let
  \begin{equation*}
    F_+ \colon \gerbes (H_1,H_0) \to \gerbes (G_1, G_0)
  \end{equation*}
  be the 2-functor given by sending the $(H_1,H_0)$-gerbe $\gQ$
  to its extension along the butterfly $[H_\bullet,E,G_\bullet]$.
\end{definition}
$F_+$ generalizes the functor $F^0_{+}$ (see~\eqref{eq:30}), and
reduces to it when $F$ arises from a strict morphism of crossed
modules.  However, note that while for a strict morphism
$f_\bullet\colon H_\bullet\to G_\bullet$ the resulting functor
$F^0_+$ reviewed in section~\ref{sec:Extension-stric} is such
that there always is an $f_\bullet$-morphism $\gQ\to F^0_+
(\gQ)$, it is not so in the current more general situation,
unless one reverts to a torsor picture.

\subsection{Induced map on non-abelian cohomology}
\label{sec:map-non-ab}

We now consider the effect of $F_+$ on cohomology.  To this end,
consider the cohomology class determined by the $(H_1,H_0)$-gerbe
$\gQ$, and let $(y,h)$ be a representative 1-cocycle with values
in $H_\bullet$, relative to a hypercover $U_\bullet \to *$.  The
class of $\gP = F_+ (\gQ)$ is obtained by applying the procedure
of section~\ref{sec:push-cohom-class} to the class of $\gQ$. More
precisely, we have:
\begin{proposition}
  \label{prop:5}
  The lift of $(y,h)$ along the butterfly, as described in
  sect.~\ref{sec:lift-cocycle-along}, provides a representative
  for the cohomology class of the $(G_1,G_0)$-gerbe $\gP$
  constructed in Theorem~\ref{thm:2}.
\end{proposition}
\begin{proof}
  The cocycle $(y,h)$ is determined by the choice of an object
  $z\in \Ob\gQ_{U_0}$, a trivialization $s$ of the $H_0$-torsor
  $\nu (z)$, and the choice of an appropriate morphism $a\colon
  d_0^*z\to d_1^*z$ over $U_1$, see the proof of
  Proposition~\ref{prop:3}.

  To prove the proposition, we show the lift of $(y,h)$ along the
  butterfly comes from a labeling of the $(H_1\times
  G_1,E)$-gerbe $\gQ'$ provided by a pair $(z', a')$, where $z'$
  is an object, and $a'\colon d_0^*z' \to d_1^*z'$ a morphism,
  respectively mapping to $z$ and $a$ under the projection $\gQ'
  \to \gQ$.  (The pair $(z',a')$ determines a non-abelian
  1-cocycle with values in $H_1\times G_1\to E$ for the gerbe
  $\gQ'$.)

  Only the construction of $z'$ and $a'$ will be carried out,
  leaving the details of the calculation that this indeed yields
  the lift of $(y,h)$ along the butterfly to the reader. In the
  process, the hypercover $U_\bullet$ will need replacing with a
  finer one, say $U'_\bullet$, by a process we have already met
  several times, now, and it will be silently done without
  further mentioning.  The need for some construction to hold
  ``locally'' will signify the need for said replacement.

  The object $z'$ can be found as follows: if $\nu\colon \gQ\to
  \tors (H_0)$ is the functor which is part of the
  $(H_1,H_0)$-gerbe structure of $\gQ$, choose a (local) lift of
  the $H_0$-torsor $\nu (z)$ to an $E$-torsor $P$, so that there
  is a $\pi$-morphism of torsors
  \begin{equation}
    \label{eq:35}
    \sigma \colon P \lto \nu (z),
  \end{equation}
  where $\pi\colon E\to H_0$. Then set $z'=(z,f, P)$, where $f$
  is the inverse of the morphism induced by $\sigma$:
  \begin{align*}
    \bar\sigma  \colon P\cprod{E} H_0 & \lto \nu (z)\\
    (p,y) &\longmapsto \sigma (p)\,y.
  \end{align*}
  
  A morphism $a' \colon d_0^*z' \to d_1^*z'$ mapping to $a\colon
  d_0^*z\to d_1^*z$ under the projection $\gQ' \to \gQ$ is of the
  form $a' = (a, \alpha)$, where $\alpha \colon d_0^*P \lto
  d_1^*P$. In fact $\alpha$ can be constructed as a (local) lift
  of $\nu (a)$ with respect to the $\pi$-morphism~\eqref{eq:35},
  so that we have a commutative diagram
  \begin{equation}
    \label{eq:36}
    \vcenter{%
      \xymatrix{%
        d_0^*P \ar[r]^{\alpha} \ar[d]_{d_0^*\sigma} & d_1^*P
        \ar[d]^{d_1^*\sigma} \\
        d_0^* \nu (z) \ar[r]_{\nu (a)} & d_1^* \nu (z)
      }
    }
  \end{equation}
  as follows. Choose $\tilde s$ of $P$ such that $\sigma (\tilde
  s)=s$, again changing $U_\bullet$ if necessary.  Indeed, note
  that the ``fiber'' $P_s=\sigma^{-1}(s)$ is a $G_1$-torsor, so
  finding $\tilde s$ amounts to a trivialization of $P_s$.  Let
  $e\in E (U_1)$ be a local lift of $y\in H_0(U_1)$ and define
  $\alpha$ as:
  \begin{equation*}
    \alpha (d_0^* \tilde s) = (d_1^*\tilde s) \, e.
  \end{equation*}
  Since $y$ is determined by the relation $\nu (a) (d_0^*s) =
  (d_1^*) y$, it is clear that $\alpha$ so defined
  satisfies~\eqref{eq:36}.
  
  Now, a further pull-back to $U_2$ determines an automorphism
  $\eta'$ of $(d_0d_1)^*z'$ such that
  \begin{equation}
    \label{eq:37}
    d_1^* a' = d_2^*a' \circ d_0^*a'\circ \eta'
  \end{equation}
  via the analog of diagram~\eqref{eq:23} in the proof of
  Proposition~\ref{prop:3}.  By construction, the projection
  $\gQ'\to \gQ$ maps $\eta'$ to the automorphism $\eta$ of
  $(d_0d_1)^*z$ obtained in the same way from $a\colon d_0^*z\to
  d_1^*z$.  It follows that $\eta' = (\eta, \epsilon)$, where
  $\epsilon$ is an automorphism of $(d_0d_1)^*P$ covering $\nu
  (\eta)$.  By using~\eqref{eq:34} we have that
  \begin{equation*}
    \shAut ((d_0d_1)^*z' \isoto (d_0d_1)^*P \cprod{E}\,
    ( H_1\times_{H_0} E),
  \end{equation*}
  so that, relative to the chosen a trivialization $\tilde s$ of
  $P$ (suitably pulled back to $U_2$), $\eta'$ is identified with
  an element of $H_1\times_{H_0} E$.  In particular, $\epsilon$
  is identified with the $E$-factor, call this particular element
  $e'\in E(U_2)$, whereas the $H_1$ factor is $h\in H_1 (U_2)$,
  which corresponds to $\eta$ via the chosen trivialization $s$
  of $\nu (z)$.  So, explicitly, the pair $(h,e')$ satisfies
  $\del (h) = \pi(e')$.  Finally, the isomorphism
  $H_1\times_{H_0}E \iso H_1 \times G_1$, identifies $(h,e')$
  with $(h,g)$, for a suitable $g\in G_1(U_2)$, or put it
  differently, $e = \kappa (h)\,\imath (g)$.

  Calculating the relation~\eqref{eq:37} with respect to the
  chosen trivializations $s$ and $\tilde s$, we find that $e$,
  $h$, and $g$ satisfy
  \begin{equation*}
    d_1^*e = d_2^*e\, d_0^*e\, \kappa (h)\, \imath (g),
  \end{equation*}
  which is the same as~\eqref{eq:16}.  Moreover, from the second
  relation of~\eqref{eq:24} applied to the pair $(a',\eta')$, or
  alternatively performing the calculation suggested at the end
  of~\ref{sec:lift-cocycle-along}, it follows that $e$, $h$, and
  $g$ also satisfy~\eqref{eq:17}, and so the 1-cocycle $(x,g)$,
  where $x=\jmath (e)$, is the lift of $(y,h)$ along the
  butterfly, as wanted.

  To complete the proof, we must make sure $(x,g)$ indeed is the
  1-cocycle arising from a labeling of the gerbe $\gP$, obtained
  from $\gQ'$ via the strict morphism $E_\bullet \to
  G_\bullet$. This is clear, since from
  section~\ref{sec:Extension-stric} we have that $\gP$ has
  locally the same objects as $\gQ'$, the functor $\mu\colon
  \gP\to \tors (G_0)$ is locally the composition of $\nu'$ with
  $\jmath_*\colon \tors (E)\to \tors (G_0)$, and the automorphism
  group of an object is locally isomorphic to $G_1$ via
  \begin{equation*}
    H_1\times_{H_0} E \iso H_1\times G_1 \lto G_1.
  \end{equation*}
\end{proof}
It follows from the previous proposition and from the arguments
in section~\ref{sec:push-cohom-class} that the class gerbe $\gP$
is therefore the image of that of $\gQ$ under $F$.  The following
is an immediate consequence of the previous results.
\begin{theorem}
  \label{thm:3}
  The gerbe $\gP$ constructed in Theorem~\ref{thm:2} is
  equivalent to $F_{**}(\gQ)$.  The two 2-functors $F_{**}$ and
  $F_{+}$ are equivalent.
\end{theorem}

\section{Commutativity conditions}
\label{sec:comm-cond}

The group law of a gr-stack may be equipped with commutativity
constraints.  Cohomology with values in such a gr-stack will inherit
corresponding structures, actually in a more rigid form due to the
process of modding out by the relation generated by (functorial)
equivalence.  Butterflies help one to obtain explicit forms for these
structures.  (Commutativity conditions for gr-stacks are thoroughly
discussed \cite{MR95m:18006,MR1702420}, see also the discussion in
\citepalias[\S 7]{ButterfliesI}.)

\subsection{Commutativity conditions and butterflies}
\label{sec:comm-cond-butt}

The very first commutativity condition one may impose on a
gr-stack is that the group law\footnote{We are going to use a
  plain symbol $m$ to denote the monoidal structure of $\grg$, in
  place of the forbidding $\otimes_\grg$ used in
  \citepalias{ButterfliesI}.}
\begin{equation}
  \label{eq:38}
  m \colon \grg\times \grg \lto \grg
\end{equation}
be \emph{braided,} that is that there be a functorial isomorphism
\begin{equation*}
  s_{x,y} \colon x \, y \lto y\, x
\end{equation*}
for each pair of objects $x,y$ of $\grg$.  Following the
convention adopted in \citepalias{ButterfliesI} (which is not the
same as refs.\ \cite{MR95m:18006,MR1702420}) we say that the
braiding is \emph{symmetric} if for all pairs of objects $x,y$ of
$\grg$ the additional condition
\begin{equation*}
  s_{y,x}\circ s_{x,y} = \id_{x\, y}
\end{equation*}
holds. In addition the symmetric braiding is \emph{Picard} if it
satisfies
\begin{equation*}
  s_{x,x} = \id_{x\, x}
\end{equation*}
for each object $x$.  A braiding is equivalent to the group law
being a \emph{morphism of gr-stacks,} rather than just a morphism
of the underlying stacks, which is the categorical analogue of
the very well-known fact that a group is abelian if and only if
its multiplication map is a group homomorphism. Therefore there
is a butterfly
\begin{equation}
  \label{eq:39}
  \vcenter{%
    \xymatrix@R-0.5em{%
      G_1\times G_1\ar[dd]_{\del\times \del}
      \ar@/_0.1pc/[dr]^\alpha  & &
      G_1 \ar@/^0.1pc/[dl]_\beta \ar[dd]^\del \\
      & P\ar@/_0.1pc/[dl]_\rho \ar@/^0.1pc/[dr]^\sigma &  \\
      G_0\times G_0 & & G_0
    }}
\end{equation}
representing the morphism~\eqref{eq:38}, see
\citepalias[7.1.3]{ButterfliesI}, once an equivalence $\grg\iso
[G_1\to G_0]\sptilde$ has been chosen.  This particular butterfly
has certain additional properties, in particular it is always
\emph{strong,} namely it always possesses a global set-theoretic
section $\tau$ of the epimorphism $\rho\colon P\to G_0\times
G_0$, so that a classical braiding map (\cite{MR1250465})
\begin{equation*}
  c\colon G_0\times G_0\lto G_1
\end{equation*}
can be obtained, see \citepalias[\S 7.1]{ButterfliesI}. The group
law of $P$ can then be described explicitly in terms of the
set-theoretic isomorphism $P\isoto G_0\times G_0\times G_1$
determined by $\tau$ and the braiding.

Depending on whether the braiding is symmetric or Picard, the
butterfly~\eqref{eq:39} satisfies extra symmetry conditions,
described in detail in \citepalias[\S 7]{ButterfliesI}.  Briefly,
if $\grg$ is braided symmetric the corresponding
butterfly~\eqref{eq:39} has the property that its pull-back under
the map that swaps the two factors in $G_\bullet \times
G_\bullet$ is isomorphic to $P$.  If in addition $\grg$ is
Picard, then the pull-back of this isomorphism to the diagonal is
the identity.

\subsection{The monoidal 2-stack of $\grg$-torsors}
\label{sec:monoidal-2-stack}

Let $\grg$ be at least braided. Since the monoidal structure of
$\grg$ is a morphism of gr-stacks, we obtain a 2-functor:
\begin{equation}
  \label{eq:40}
  m_*\colon \tors (\grg) \times \tors (\grg)  \lto \tors (\grg) 
\end{equation}
where we have used the identification $\tors (\grg \times \grg)
\iso \tors (\grg) \times \tors (\grg)$.  Thus, $m_*$ assigns to
the $\grg \times \grg$-torsor $(\stX,\stX')$ the $\grg$-torsor
$(\stX,\stX')\cprod{\grg\times\grg} \grg$.

By the theory of section~\ref{sec:extension-butterfly} the gerbe
counterpart of~\eqref{eq:40} is the 2-functor
\begin{equation}
  \label{eq:41}
  m_+\colon \gerbes (G_1,G_0) \times \gerbes (G_1,G_0)  \lto
  \gerbes (G_1,G_0) 
\end{equation}
given by the lift of the gerbe $(\gP,\gP')$ along the
butterfly~\eqref{eq:39}.

A full investigation of the monoidal structure~\eqref{eq:40}
or~\eqref{eq:41} is beyond the scope of the present work, but it
is necessary to at least point out that it is the entire
collection (in this case: 2-gerbe) of geometric objects itself
that acquires a (weak) group structure. The one on cohomology is
then obtained by considering equivalence classes, and it is
examined in the next section.

\subsection{Group structures on cohomology and butterflies}
\label{sec:group-struct-cohom}

If $\grg$ is at least braided, its monoidal
structure~\eqref{eq:38} induces morphisms
\begin{equation}
  \label{eq:42}
  m_*\colon 
  \H^i(*,\grg) \times \H^i(*,\grg) \lto \H^i(*,\grg),
\end{equation}
by the mechanisms expounded both in \citepalias{ButterfliesI}
(for degree $i \leq 0$) and in the present work (for degree
$i=0,1$).  The morphism~\eqref{eq:42} is obtained starting from
either~\eqref{eq:40} or~\eqref{eq:41} and using functoriality.

At the level of representing cocycles, the group laws~\eqref{eq:42}
can be computed by applying the lifting along the
butterfly~\eqref{eq:39} described in
section~\ref{sec:lift-cocycle-along} (By the observation in
remark~\ref{rem:4}, it applies equally well to 0-cocycles, \ie descent
data for objects of gr-stacks).  The weak form of the group law for
$\grg$ translates into a standard rigid one for the $m_*$, including
the case $i=1$.

We collect the main facts in the following
\begin{proposition}
  \label{prop:6}
  Let $\grg$ be a braided gr-stack.
  \begin{enumerate}
  \item \label{item:3} $\H^0(*,\grg)$ is an abelian group;
  \item \label{item:4} $\H^1(*,\grg)$ is a group;
  \item \label{item:5} If in addition $\grg$ is symmetric,
    $\H^1(*,\grg)$ is an abelian group.
  \end{enumerate}
\end{proposition}
\begin{proof}[Sketch of the proof]
  The result is quite well-known, so we only sketch the main
  ideas.
  
  For \ref{item:3}, given that $\H^0(*,\grg) \iso \pi_0 (\grg
  (*))$, the result is obvious (it follows immediately from the
  weak group law of $\grg$).  As noted, for case~\ref{item:4},
  that is $\H^1(*,\grg)$, it follows from either morphism in
  section~\ref{sec:monoidal-2-stack} and functoriality.

  More interesting is the case of a symmetric gr-stack.  It was
  proved in \citepalias[Propositions 7.2.2 and
  7.2.3]{ButterfliesI} that the symmetry condition is equivalent
  to the braiding being a \emph{2-morphism}
  \begin{equation*}
    s \colon m \Longrightarrow
    m \circ T \colon \grg\times \grg \lto \grg
  \end{equation*}
  of gr-stacks, where $T$ is the swap functor.  Passing to
  cohomology and using~\eqref{eq:42} yields the commutative
  structure
  \begin{equation*}
    \xymatrix{%
      \H^1(*,\grg) \times \H^1(*,\grg) \ar[r]^(.6)\sim \ar[d]_{T_*} &
      \H^1(*,\grg\times \grg)  \ar[d]_{T_*}
      \ar[r]^(.6){m} & 
      \H^1(*,\grg) \ar@{=} [d]\\
      \H^1(*,\grg) \times \H^1(*,\grg) \ar[r]^(.6)\sim &
      \H^1(*,\grg\times \grg)  \ar[r]^(.6){m} &
      \H^1(*,\grg) 
    }
  \end{equation*}\qedhere
\end{proof}

\subsection{Explicit cocycles}
\label{sec:explicit-cocycles}

Besides ``explaining'' how the first non-abelian cohomology group with
values in a crossed module acquires a group structure, with the
butterfly we can calculate explicit formulas for the product.  The
computations involved are tedious and straightforward overall, so we
will not dwell on the details and only report the main formulas.

As already observed the butterfly~\eqref{eq:39} is strong, so the
group law of $P$ can be explicitly described in terms of the
set-theoretic isomorphism $P\iso G_0\times G_0\times G_1$ and the
braiding $c$ as
\begin{equation*}
  (x_0,y_0,g_0)\, (x_1,y_1,g_1) = (x_0x_1, y_0y_1,
  c(x_1,y_0)^{y_1}g_0^{y_0y_1}g_1),
\end{equation*}
with $x_0,x_1,y_0,y_1\in G_0$, and $g_0,g_1\in G_1$.  In the
foregoing the strong set-theoretic section $\tau\colon G_0\times
G_0 \to P$ is obviously of the form
\begin{equation*}
  \tau (x,y) = (x,y,1),
\end{equation*}
with $x,y\in G_0$.  In fact, all the maps in~\eqref{eq:39} have
explicit descriptions in these coordinates, and their form will
be left as an exercise to the interested reader; here we only
mention that $\sigma\colon P\to G_0$ has the form
\begin{equation*}
  \sigma (x,y,g) = x\, y\, \del g.
\end{equation*}
Note that the composition with $\tau$ gives the multiplication
map of $G_0$, which is of course not a homomorphism.\footnote{In
  this way one arrives at the standard interpretation of the
  braiding map as the isomorphism relating the multiplication map
  and its swapped version.}  The two main computations are as
follows.

\subsubsection*{Degree zero}

Assume two global objects $X,X'\in \Ob \grg (*)$ are represented
by zero-cocycles (descent data) $(x,g)$ and $(x',g')$ relative to
some common (hyper)cover $U_\bullet \to *$.  Here $x,x' \in G_0
(U_0)$ and $g,g'\in G_1 (U_1)$.  The object $(X,X')$ of
$\grg\times \grg$ is represented by the direct product of the
corresponding cocycles.  Applying the procedure of
section~\ref{sec:lift-cocycle-along} (adapted to 0-cocycles, as
per Remark~\ref{rem:4}) one finds that the image of $(X,X')$
under the multiplication map~\eqref{eq:42} is represented by the
cocycle
\begin{equation*}
  (xx', g^{d_1^*x}g').
\end{equation*}
This formula coincides with the one for the group law of the
gr-stack $\grg$ expressed in terms of descent data found in
\citepalias[3.4.3]{ButterfliesI}.  So the lift along the
butterfly computes exactly the same (abelian) group law as
induced by the braided structure on $\grg$.
\begin{remark}
  A priori there appear to be \emph{two} group laws on
  $\H^0(*,\grg)$.  One inherited from the monoidal structure of
  $\grg$, while the second is $m_*$ in~\eqref{eq:42}.  One is a
  homomorphism of the other, so by the classical argument they
  coincide, and the resulting structure is abelian.
\end{remark}

\subsubsection*{Degree one}

Assume now $\gP,\gP'$ are two gerbes bound by the crossed module
$G_\bullet$.  Recycling symbols, assume they are represented by
1-cocycles $(x,g)$ and $(x',g')$ relative to some common
(hyper)cover $U_\bullet \to *$.  This time $x,x' \in G_0 (U_1)$
and $g,g'\in G_1 (U_2)$.  The product gerbe $\gP\times \gP'$ is
represented by the direct product of the corresponding cocycles.
Applying again the procedure of
section~\ref{sec:lift-cocycle-along} the gerbe $m_{+} (\gP\times
\gP')$ of section~\ref{sec:extension-butterfly} (see in
particular Definition~\ref{def:6}) is represented by a 1-cocycle
relative to $U_\bullet$ given by the expression:
\begin{equation}
  \label{eq:43}
  \bigl( x \, x' , 
  c(d_0^*x,d_2^*x')^{-d_0^*x'} g^{d_2^*x'\, d_0^*x'}g' \bigr).
\end{equation}
We could have used $\grg$-torsors $\stX$ and $\stX'$ to arrive at
the same conclusion. In particular, if $(x,g)$ and $(x',g')$ are
assumed to be 1-cocycles corresponding to $\stX$ and $\stX'$,
then the 1-cocycle of expression~\eqref{eq:43} represents the
$\grg$-torsor $(\stX \times \stX')\cprod{\grg\times \grg} \grg$.

In summary, modulo the appropriate notion of equivalence,
expression~\eqref{eq:43} gives an explicit form to the group
law~\eqref{eq:42} when $i=1$.

If $\grg$ is braided symmetric, the geometric condition on the
butterfly~\eqref{eq:39} translates into the standard notion that
the braiding map satisfies the symmetry condition
$c(x,y)=c(y,x)^{-1}$. In this instance it is possible to
explicitly verify that $\H^1(*,\grg)$ becomes an abelian group;
exchanging the role of $(x,g)$ and $(x',g')$ in
expression~\eqref{eq:43} leads to a 1-cocycle which can be seen
to be equivalent to the original one.  We omit the details.

\section{Butterflies and extensions}
\label{sec:butt-extens}

Group extensions and non-abelian cohomology in degree one have a
close relationship, which one can trace from Dedeker's classical
approach based on cocycle calculations, to Grothendieck's and
Breen's more geometric one, where the category of extensions
\begin{equation*}
  1\lto G\lto E\lto \Gamma \lto 1
\end{equation*}
of the topos $\T$ is given geometric meaning by showing its
equivalence to that a morphism of gr-stacks
\begin{equation*}
  \Gamma \lto \bitors (G).
\end{equation*}
$\bitors (G)$ is the gr-stack associated to the crossed module
\begin{math}
  G \to \Aut (G),
\end{math}
and $\Gamma$ is considered as a gr-stack in the obvious way.
These ideas fit very well within the butterfly framework.

\subsection{The Schreier-Grothendieck-Breen theory of extensions}
\label{sec:schr-groth-breen}

Following ref.\ \cite[\S 8.11]{MR92m:18019}, consider an
\emph{extension of $\Gamma$ by the crossed module $G_1\to G_0$,}
a notion due to Dedecker and defined by the following commutative
diagram:
\begin{equation}
  \label{eq:44}
  \vcenter{%
    \xymatrix{%
      1 \ar[r] & G_1 \ar[r]^\imath \ar[d]_\del &
      E \ar[r]^\pi \ar@/^/[dl]^\jmath & \Gamma \ar[r] & 1 \\
      & G_0
    }
  }
\end{equation}
where the map $\jmath\colon E\to G_0$ is subject to the
additional condition
\begin{equation}
  \label{eq:45}
  e^{-1}\imath (g)e = \imath (g^{\jmath (e)}).
\end{equation}
We recognize~\eqref{eq:45} as the first relation in~\eqref{eq:7},
as well as \cite[equation (8.11.2)]{MR92m:18019}, after the
obvious changes due to the different conventions adopted in this
paper.

The trivial extension corresponds to $E=\Gamma \ltimes G_1$,
where $\Gamma$ acts on $G_1$ via a homomorphism $\xi\colon
\Gamma\to G_0$ and the action of $G_0$ on $G_1$, whereas $\jmath$
is given set-theoretically as
\begin{equation*}
  \jmath (x,g) = \xi (x) \, \del g,
\end{equation*}
for $x\in\Gamma$ and $g\in G_1$.

A comparison with diagram~\eqref{eq:6} suggests
diagram~\eqref{eq:44} ought to be considered as a ``one-winged
butterfly,'' namely a butterfly diagram from the crossed module
$[1\to \Gamma]$ to $[G_1\to G_0]$.  Therefore, by the results in
\citepalias[\S 4 and \S 5]{ButterfliesI}, recalled in
section~\ref{sec:butt-weak-morph}, the extension~\eqref{eq:44}
corresponds to a morphism of gr-stacks
\begin{equation*}
  F_E\colon \Gamma \lto \grg
\end{equation*}
where $\grg\iso [G_1\to G_0]\sptilde$.  The form of this morphism
is as follows.  If $x\colon U\to \Gamma$ is a point, it follows
from \citepalias[\S 4.3]{ButterfliesI} (see also
section~\ref{sec:computing-f_colon-h1} for a quick review), that
it maps to the $(G_1\rvert_U,G_0\rvert_U)$-torsor
\begin{equation*}
  \shHom_1(1,E)_x \iso x^*E\coin E_x.
\end{equation*}
This retrieves the expression \citep[8.2.2]{MR92m:18019}.
Observe also that \eqref{eq:45} is none other than the expression
of the left $G_1$-action on $x^*E$ in terms of the right one
(cf.\ section~\ref{sec:crossed-modules-gr}).  In this language a
trivial extension corresponds to a split butterfly.  Note also
that for a split extension the $(G_1\rvert_U,G_0\rvert_U)$-torsor
$x^*E$ is isomorphic to $(G_1\rvert_U,x)$.

The obvious notion of morphism of extensions of the
form~\eqref{eq:44} is clearly the same as that of morphism of
one-winged butterflies, in other words an isomorphism $\phi\colon
E\to E'$ of group objects compatible with~\eqref{eq:44}.  With
reference to the notation used elsewhere in this series (see,
e.g. section~\ref{sec:butt-weak-morph}) we have
\begin{equation*}
  \catExt(\Gamma, G_1\to G_0) \coin
  \cat{B} (\Gamma,G_\bullet),
\end{equation*}
where the left-hand side denotes the category (in fact, the
groupoid) of extensions of the form~\eqref{eq:44}, and the
right-hand side the one of butterflies.  It immediately follows
from Theorem~\ref{thm:1} that there is an equivalence of
categories
\begin{equation}
  \label{eq:46}
  \catExt(\Gamma,G_1\to G_0) \lisoto \catHom(\Gamma,\grg).
\end{equation}
There is also the fibered analog of the preceding construction.
Again from \citepalias[\S 4 and \S 5]{ButterfliesI} (see also the
summary in section~\ref{sec:bicat-cross-modul}), and using the
same notation, we obtain the following analog of \cite[Lemme
8.3]{MR92m:18019}:
\begin{lemma}
  \label{lem:7}
  There is an equivalence
  \begin{equation*}
    \shcatExt(\Gamma,G_1\to G_0) \lisoto \shcatHom(\Gamma,\grg),
  \end{equation*}
  where the left-hand side is the stack whose fiber over $U$ is
  $\catExt (\Gamma\rvert_U, G_\bullet\rvert_U)$.
\end{lemma}
The cohomological classification of the extensions is obtained by
applying $\pi_0$ to~\eqref{eq:46},
\begin{equation*}
  \Ext(\Gamma,G_1\to G_0) \lisoto \Hom(\Gamma,\grg),
\end{equation*}
and rephrasing the right-hand side in terms of the non-abelian
cohomology of the classifying object $\B\!\Gamma$.  Briefly, the
group structure of $\Gamma$ is encoded by diagram 8.1.2 of
\citep{MR92m:18019}, which we write in the form
\begin{equation}
  \label{eq:47}
  \gamma \colon d_1^*E \lisoto d_2^*E\, \cprod{G_1}\, d_0^*E,
\end{equation}
subject to the coherence condition for $\gamma$ expressing the
associativity of the group law.  Pulling back by $x\colon U\to
\Gamma$, and then $d_0^*x$, $d_1^*x$, $d_1^*x$, we can
see~\eqref{eq:47} plus the coherence condition for $\gamma$
define a 1-cocycle on $\B\!\Gamma$ with values in $\grg$.  By a
reasoning entirely analogous to the one of
section~\ref{sec:computing-f_colon-h1}, we can compute the class
with values in the crossed module $G_\bullet$, thereby obtaining
the sought-after element in $\H^1(\B\!\Gamma, \grg)$. Thus we
have:
\begin{proposition}[{\citealp[Proposition 8.2]{MR92m:18019}}]
  \label{prop:7}
  There is a functorial isomorphism of sets
  \begin{equation*}
    \Ext(\Gamma,G_1\to G_0) \lisoto \H^1(\B\!\Gamma, \grg).
  \end{equation*}
\end{proposition}
Functoriality is built into the butterfly representation of
morphisms of gr-stacks.

\subsection{Remarks on extensions by commutative crossed modules}
\label{sec:remarks-extens-comm}

We can combine the idea of extension by a crossed
module~\eqref{eq:44} with the conditions studied in
section~\ref{sec:comm-cond}.  In this situation the first
non-abelian cohomology set
\begin{math}
  \H^1(\B\!\Gamma, \grg)
\end{math}
acquires a group structure, possibly abelian if $G_\bullet$ is
symmetric or Picard.

\subsubsection*{Baer sums}

The explicit cocycle multiplication formula~\eqref{eq:43} could easily
be translated in terms of group cohomology.  This is easier in the
case of a strong butterfly, that is for an extension~\eqref{eq:44}
possessing a global set-theoretic section $s\colon \Gamma\to E$, and
it is left as an exercise to the reader.

There is a more interesting ``butterfly explanation'' of the
existence of the product; while the basic mechanism is the one
already explained in section~\ref{sec:comm-cond}, the translation
in terms of group cohomology gives it a slightly different flavor
that further underscores the role of butterfly diagrams.  The
procedure outlined below is the analog in the context of
non-abelian cohomology of the standard Baer sum of extensions in
ordinary homological algebra (see \cite{maclane:hom}).

From two extensions of type~\eqref{eq:44}, we can form the direct
product (drawn with a different orientation) one-winged
butterfly:
\begin{equation}
  \label{eq:48}
  \vcenter{%
    \xymatrix@R-0.5em{%
      & &  G_1\times G_1 \ar@/^0.1pc/[dl]_{(\imath,\imath')}
      \ar[dd]^{(\del,\del)}\\
      & E\times E' \ar@/_0.1pc/[dl]_{(\pi,\pi')}
      \ar@/^0.1pc/[dr]^{(\jmath,\jmath')} &  \\
      \Gamma\times \Gamma & & G_0\times G_0
    }}
\end{equation}
which then can be composed with~\eqref{eq:39}, which encodes the
monoidal structure, to yield
\begin{equation*}
  \xymatrix@R-0.5em{%
    & &  G_1\times G_1 \ar@/^0.1pc/[dl]_{(\imath,\imath)}
    \ar[dd]^{(\del,\del)} \ar@/_0.1pc/[dr]^\alpha & &
    G_1 \ar@/^0.1pc/[dl]_\beta \ar[dd]^\del \\
    & E\times E' \ar@/_0.1pc/[dl]_{(\pi,\pi)}
    \ar@/^0.1pc/[dr]^{(\jmath,\jmath)} & &
    P\ar@/_0.1pc/[dl]_\rho \ar@/^0.1pc/[dr]^\sigma &  \\
    \Gamma\times \Gamma & & G_0\times G_0 & & G_0
  }
\end{equation*}
that is, according to \citepalias[\S 5.1]{ButterfliesI},
\begin{equation*}
  \xymatrix@R-0.5em{%
    & & G_1 \ar@/^0.1pc/[dl] \ar[dd]^\del \\
    &{\displaystyle \bigl( E\times E' \bigr) \times_{G_0\times G_0}^{G_1\times
        G_1} P }
    \ar@/_0.1pc/[dl]
    \ar@/^0.1pc/[dr] &  \\
    \Gamma\times \Gamma & & G_0
  }
\end{equation*}
which is then pulled back to $\Gamma$ via the diagonal
homomorphism $\Delta\colon \Gamma\to \Gamma\times \Gamma$.  The
overall picture for the product is as follows:
\begin{equation*}
  \xymatrix@R-0.5em{%
    1 \ar[rr] \ar[dd] & & 1 \ar[dd] \ar@/_0.1pc/[dr] & &
    G_1 \ar@/^0.1pc/[dl] \ar[dd]^\del \\
    & & &{\displaystyle \bigl( E\times E' \bigr) \times_{G_0\times G_0}^{G_1\times
        G_1} P }
    \ar@/_0.1pc/[dl]
    \ar@/^0.1pc/[dr] &  \\
    \Gamma \ar[rr]^\Delta & & \Gamma\times \Gamma & & G_0
  }
\end{equation*}
The composition expressed by the above diagram is the full
butterfly diagram expressing the product structure on the first
cohomology with coefficients in $\grg$.  Thus we obtain a
monoidal structure on the category $\catExt(\Gamma,G_\bullet)$.

\subsubsection*{Abelian structure on $\H^1$}

If $\grg$ (or equivalently $G_\bullet$) is symmetric, the
butterfly~\eqref{eq:39} is isomorphic to itself under pull-back
by the morphism $T$ that switches the factors. By
\citepalias[\S 7.2.4]{ButterfliesI} this means there exists $\psi
\colon P\isoto P$ such that:
\begin{equation*}
  \xymatrix{%
    P \ar[d] \ar[r]^\psi & P \ar[d] \\
    G_0\times G_0 \ar[r]^T & G_0\times G_0
  }
\end{equation*}
compatible with all the morphisms in~\eqref{eq:39}.  The same
kind of swap of course exchanges the factors in the
butterfly~\eqref{eq:48}.  Therefore there is a diagram of
juxtaposed butterflies
\begin{equation*}
  \xymatrix{%
    \Gamma \ar@{=}[d] \ar[r]^\Delta &
    \Gamma\times \Gamma \ar[d]^T &
    E\times E' \ar[l] \ar[d]^T \ar[r] &
    G_\bullet \times G_\bullet \ar[d]^T &
    P \ar[l] \ar[d]^\psi \ar[r] &
    G_\bullet \ar@{=}[d] \\
    \Gamma \ar[r]^\Delta &
    \Gamma\times \Gamma  &
    E'\times E \ar[l] \ar[r] &
    G_\bullet \times G_\bullet &
    P \ar[l] \ar[r] &
    G_\bullet  }
\end{equation*}
which leads to a morphism of one-winged butterflies
\begin{equation*}
  \xymatrix@R-1pc{%
    & \Delta^*\bigl( (E\times E') \times_{G_0\times G_0}^{G_1\times
      G_1} P \bigr)  \ar@/_0.2pc/[dl] \ar[dd] \ar@/^0.2pc/[dr] \\
    \Gamma && G_\bullet \\
    & \Delta^*\bigl( (E'\times E) \times_{G_0\times G_0}^{G_1\times
      G_1} P \bigr)  \ar@/^0.2pc/[ul] \ar@/_0.2pc/[ur] 
  }
\end{equation*}
from $\Gamma$ to $G_\bullet$.  This provides a purely
diagrammatic proof that the group structure of
$\H^1(\B\!\Gamma,\grg)$ is abelian when $\grg$ is symmetric.  At
the level of diagrams, it is a braiding on the category
$\catExt(\Gamma,G_\bullet)$.

\subsection{Butterflies, extensions, and simplicial morphisms}
\label{sec:butt-simpl-morph}

Consider again a generic morphism $F\colon \grh\to \grg$ of
gr-stacks and the corresponding butterfly~\eqref{eq:6}.  Using a
sheafified nerve construction, $F$ corresponds to a simplicial
map
\begin{equation}
  \label{eq:49}
  \W\smp{H}_\bullet \lto \W\smp{G}_\bullet,
\end{equation}
via the map $\smp{H}_\bullet \to \smp{G}_\bullet$ in the sense of
$A_\infty$-spaces, thanks to considerations analogous to those of
\citep[\S 8.5]{MR92m:18019}.  In the set-theoretic case this
simplicial map is the starting point for the definition of
weak-morphism of crossed module, which is then \emph{computed} by
a butterfly diagram. In the sheaf-theoretic context the starting
point for the definition of weak morphism is different (See the
discussion in \citepalias[\S 4.2]{ButterfliesI}).  Thus, it is of
some interest to re-obtain the simplicial map in the present
context.

Rather than appealing to $A_\infty$-geometry, we sketch a
different way to arrive at the same conclusion, as follows.  If
in the butterfly~\eqref{eq:6} we isolate the ``one-winged'' one,
\begin{equation}
  \label{eq:50}
  \vcenter{%
    \xymatrix@R-0.5em{%
      & &  G_1 \ar@/^0.1pc/[dl]_\imath \ar[dd]^\del\\
      & E\ar@/_0.1pc/[dl]_\pi \ar@/^0.1pc/[dr]^\jmath &  \\
      H_0 & & G_0
    }}
\end{equation}
analogous to~\eqref{eq:44}, we obtain a class in $\H^1(\B\! H_0,
\grg)$, corresponding to a well-defined morphism
\begin{equation*}
  H_0\lto \grg,
\end{equation*}
in the sense of gr-stacks.  Thus, the underlying geometric object
to the extension~\eqref{eq:50} is a $\grg$-torsor, or
equivalently, a gerbe bound by $G_\bullet$, over $\B\!H_0$.

Next, the standard pull-back (see \cite{maclane:hom}) of the
extension~\eqref{eq:50} to $H_1$ via $\del \colon H_1\to H_0$ is
trivial, due to the existence of the homomorphism $\kappa\colon
H_1\to E$ in the full butterfly~\eqref{eq:6}.  It follows that
the class of the extension~\eqref{eq:50} dies under the pull-back
map
\begin{equation}
  \label{eq:51}
  (\B\!\del)^*\colon
  \H^1(\B\! H_0, \grg) \lto \H^1(\B\! H_1, \grg).
\end{equation}
The condition that the pullback of the cocycle corresponding to
the extension~\eqref{eq:50} vanish leads to an explicit
simplicial map~\eqref{eq:49}.  The actual computation via
cocycles is uneventful and quite laborious, so we omit it.

More interesting is the geometric reason, which we record in the
following informal assertions---not all verification having being
carried out.  Essentially, the $\grg$-torsor over $\B\! H_0$
defined by the extension~\eqref{eq:50} ``descends'' to
$\W\smp{H}_\bullet$ along the map $\B\!H_0\to \W\smp{H}_\bullet$.
\begin{assertion}
  The vanishing of the image of the class of the
  extension~\eqref{eq:50} under the map~\eqref{eq:51} determines
  2-descent data for the $\grg$-torsor determined by the
  extension~\eqref{eq:50} relative to the map $\B\!H_0\to
  \W\smp{H}_\bullet$.
\end{assertion}
\begin{proof}[Sketch of the proof]
  Consider the augmented (bi)simplicial object
  \begin{equation*}
    U_{\bullet \bullet} =
    \cosk_0 \bigl( \B\!H_0\to
    \W\smp{H}_\bullet \bigr) \colon 
    \xymatrix@1{%
      \dotsm \ar@<0.9ex>[r] \ar[r] \ar@<-0.9ex>[r] & 
      \B\!H_0 \times_{\smash[b]{\W\smp{H}_\bullet}} \B\!H_0
      \ar@<0.6ex>[r] \ar@<-0.6ex>[r] & \B\!H_0 \ar[r] &
      \W\smp{H}_\bullet
    }
  \end{equation*}
  where the first index is the ``external'' one, whose face maps
  are explicitly drawn above.  We compute $B\!H_0
  \times_{\smash[b]{\W\smp{H}_\bullet}} B\!H_0 \iso \B (
  H_0\ltimes H_1)$, and so on, therefore $U_{\bullet\bullet}$ is
  equivalent to $\B$ applied degree-wise to $\smp{H}_\bullet$:
  \begin{equation*}
    \xymatrix@1{%
      \dotsm \ar@<1.4ex>[r] \ar@<0.5ex>@{}[r]|(0.2)\vdots \ar@<-1.3ex>[r]&
      \B (H_0\ltimes (H_0\ltimes H_1))
      \ar@<0.9ex>[r] \ar[r] \ar@<-0.9ex>[r] & 
      \B (H_0\ltimes H_1)
      \ar@<0.6ex>[r] \ar@<-0.6ex>[r] & B\!H_0\ar[r] &
      \W\smp{H}_\bullet
    }
  \end{equation*}
  The face maps are actually induced by those of
  $\smp{H}_\bullet$.  Note that the diagonal of the above
  bisimplicial object is equivalent to $\W\!\smp{H}_\bullet$.

  The extension~\eqref{eq:50} determines a bitorsor cocycle of
  the type~\eqref{eq:47} which we write as:
  \begin{equation*}
    \gamma_{x,y}\colon E_{xy} \lisoto E_x \cprod{G_1} E_y,
  \end{equation*}
  for points $x,y$ of $H_0$.  The class of this cocycle is trivial
  under the pull-back~\eqref{eq:51}, and moreover we know the
  pulled-back extension is actually a \emph{direct} product, rather
  than merely a semi-direct one, since the composition $\jmath\circ
  \kappa$ is trivial in the full butterfly.  A moment's thought
  reveals the $(G_1,G_0)$-torsor determined by a direct product
  extension is in fact trivial, \ie of the form $(G_1,1)$, hence we
  must have coherent isomorphisms
  \begin{equation*}
    \delta_h \colon E_{\del h} \lisoto G_1,
  \end{equation*}
  where of course $E_{\del h}$ is the ``value'' of the pulled
  back cocycle at $h$.

  At a point $(y,h)$ of $H_0\ltimes H_1$, the pull-backs of $E$
  along the two face maps
  \begin{equation*}
    d_i\colon H_0\ltimes H_1 \lto H_0,\quad i=0,1,
  \end{equation*}
  $d_0(y,h)=y\del h$, and $d_1(y,h)= y$, are:
  \begin{equation*}
    d_0^*E_{(y,h)} = E_{y\del h}, \quad
    d_1^*E_{(y,h)} = E_y.
  \end{equation*}
  Using the cocycle condition and the triviality argument above,
  we have an isomorphism
  \begin{equation*}
    E_{y\del h} \xrightarrow{\gamma_{y,\del h}} 
    E_{y} \cprod{G_1} E_{\del h} \xrightarrow{1\wedge \delta_h} 
    E_y
  \end{equation*}
  at each point $(y,h)$ of $H_0\ltimes H_1$.  Thus, we have
  obtained an isomorphism of extensions, and hence of
  $\grg$-torsors, or again gerbes bound by $G_\bullet$, over the
  first stage $U_{1\bullet}$.

  Similar arguments, this time using the coherence of $\gamma$
  and $\delta$, would show the axioms of a 2-descent datum with
  respect to $\B\! H_0\to \W\!\smp{H}_\bullet$ are satisfied.
\end{proof}
Let us denote by $\gerbe{E}$ the descended gerbe over
$\W\!\smp{H}_\bullet$.  Finally we have:
\begin{assertion}
  The class of $\gerbe{E}$ determines the simplicial
  map~\eqref{eq:49}.
\end{assertion}
\begin{proof}[Sketch of the proof]
  After sections~\ref{sec:tors-non-ab-cohom}
  and~\ref{sec:gerbes-bound-crmod}, the class of a gerbe is
  effectively a simplicial map of the sought-after type.
\end{proof}

\bibliographystyle{halpha}
\bibliography{general}%

\end{document}